\documentclass{m2an_cleaned}
\usepackage{graphicx}
\usepackage{xcolor} 
\usepackage{amssymb,amsmath} 
\usepackage{hyperref}
\usepackage{color}
\usepackage[percent]{overpic}
\graphicspath{{./figures/}{./}}

\newcommand{\R}{\mathbb R}
\newcommand{\Dx}{\Delta x}
\newcommand{\Dt}{\Delta t}

\newcommand{\Np}{N_{\text{P}}}
\newcommand{\Nppc}{N_{\text{PC}}}
\newcommand{\Pp}{\mathcal P_{\alpha}}
\renewcommand{\P}{\mathcal P}
\newcommand{\M}{\mathcal M}
\newcommand{\GG}{\mathcal G}
\newcommand{\CFL}{\beta}

\renewcommand{\tt}{\frac{\theta}{1-\theta}}

\renewcommand{\S}{\mathcal S}
\newcommand{\W}{\mathcal W}
\newcommand{\V}{\mathcal V}

\newcommand{\ODEsolver}{\mathcal S_{\textup{ODE}}}
\newcommand{\thetauno}{\lambda}
\newcommand{\thetadue}{\mu}
\newcommand{\YES}{\Phi}
\newcommand{\Eref}{E^1_{\textup{ref}}}

\newcommand{\MMx}{1}
\newcommand{\mMsinx}{3}
\newcommand{\mMx}{2}
\newcommand{\mMlwr}{4}

\begin{document}
%%-----------------------------
%%      the top matter
%%-----------------------------
\title{Blended Numerical Schemes for the Advection Equation \\ and Conservation Laws}\thanks{This work has been partially supported by the INdAM-GNCS 2015 project ``Metodi numerici semi-impliciti e semi-lagrangiani per sistemi iperbolici di leggi di bilancio''}

\author{Simone Cacace}\address{Dipartimento di Matematica, Sapienza -- Universit\`a di Roma, Rome, Italy (cacace@mat.uniroma1.it)}

\author{Emiliano Cristiani}\address{Istituto per le Applicazioni del Calcolo ``M. Picone'', Consiglio Nazionale delle Ricerche, Rome, Italy (e.cristiani@iac.cnr.it)}

\author{Roberto Ferretti}\address{Dipartimento di Matematica e Fisica, Universit\`a Roma Tre, Rome, Italy (ferretti@mat.uniroma3.it)}

\date{May 2016}

\begin{abstract} 
In this paper we propose a method to couple two or more explicit numerical schemes approximating 
the same time-dependent PDE, aiming at creating new schemes which inherit advantages of the original ones. We consider both advection equations and nonlinear conservation laws.
By coupling a macroscopic (Eulerian) scheme with a microscopic (Lagrangian) scheme, we get a new kind of multiscale numerical method. 
\end{abstract}

\maketitle

\section{Introduction}\label{sec:intro}
In this paper we propose and analyse a method to couple two or more numerical schemes which approximate the same time-dependent PDE. The result is a new scheme which inherits advantages (and drawbacks) of the original ones. In what follows, we will discuss how the combination might allow to improve each of the schemes, and try to derive practical coupling strategies. Roughly speaking, we combine the existing schemes by a \emph{convex combination} in such a way that each scheme is influenced by the others. In the end, we get as many approximate solutions as the number of the considered schemes, but all of them are modified by the combination. In the basic version of the algorithm the coupling is performed everywhere and at any time. Nevertheless, local (in space or in time) couplings are also possible.
 
For illustrative purposes, we focus on the coupling of two schemes only. We consider both the case of two macroscopic (Eulerian) schemes and, more interesting, the case of a microscopic (Lagrangian, particle-based) scheme blended with a macroscopic scheme. In the latter case we get a new kind of multiscale numerical method. As a test problem, we consider the following advection equation in conservative form
\begin{equation}\label{AE}
\left\{
\begin{array}{ll}
u_t(x,t)+(A(x)\ u(x,t))_x=0, & x\in \R^d, \quad t>0, \\ [2mm]
u(x,0)=\bar u(x), & x\in \R^d,
\end{array}
\right .
\end{equation}
where $\bar u:\R^d\to\R$ is the initial datum and $A:\R^d\to\R^d$ is a bounded and Lipschitz continuous vector field. 
This equation describes the transport of a certain amount of mass in $\R^d$, drifted with velocity $A$. The function $u$ represents the \emph{density} of such a mass, so that the mass in any set $E\subseteq\R^d$ at any time $t$ is given by
$$
\mathcal{M}(E,t)=\int_E u(x,t)\ dx.
$$
In order to introduce a multiscale framework, it is useful to recast the problem at microscopic scale. Denoting by $\P\in\R^d$ the position of a generic mass particle, the Lagrangian counterpart of \eqref{AE} is the ODE
\begin{equation}\label{ODE1p}
\dot\P=A(\P).
\end{equation}
with suitable initial condition $\P(0)$.

The extension to the one-dimensional nonlinear conservation law
\begin{equation}\label{CL}
\left\{
\begin{array}{ll}
u_t(x,t)+f(u(x,t))_x=0, & x\in \R, \quad t>0, \\ [2mm]
u(x,0)=\bar u(x), & x\in \R,
\end{array}
\right .
\end{equation}
will be also discussed.

\medskip

The coupling idea investigated in this paper can be pursued whenever one has at his disposal two or more numerical schemes with \emph{complementary advantages and drawbacks}, resulting in a blended scheme which shows a combination of these features. The idea is mainly inspired by the multiscale technique introduced in \cite{cristiani2011MMS}, which can be applied whenever one has two equations modeling the same physical phenomenon at different scales, and the velocity field explicitly appears in both equations. The blend is obtained by a convex combination of the two velocity fields, an then both the macroscopic density and the microscopic particles are transported by means of the new common multiscale velocity field, which, in turn, depends on both density and positions of the particles. Along the same line, the paper \cite{cristiani2015JCSMD} blends the Brownian motion with the heat equation. Note that in these papers the goal is to improve the description of some physical phenomenon, by catching at large scale some effects ultimately triggered by small scale features, otherwise invisible. Conversely, in the present paper we blend directly the solutions of the schemes rather than the velocity fields, with the goal of improving the quality of the numerical approximation.

\subsubsection*{Relevant literature} 
The literature about numerical schemes for the advection equation, including multiscale numerical methods, is huge and a complete review of it is out of the scope of the paper. We will refer the interested reader to the book \cite{levequebook} for an introduction to the equations \eqref{AE} and \eqref{CL}, and the basic Eulerian numerical schemes to solve them. We also mention the book \cite{weinanbook} for a survey on the best known techniques to couple microscopic and macroscopic descriptions of physical phenomena. 
In the following we only discuss the schemes that can be put in relation with the one proposed in this paper. 
\begin{itemize}
\item The proposed scheme follows the same philosophy of Crank-Nicolson method and $\theta$-methods in general. However, there are two important differences: first, our method is explicit and combines two explicit schemes, while $\theta$-methods always combine the explicit and the implicit version of the same scheme. Second, our method has two free parameters instead of one. 

\item Filtered (hybrid) schemes \cite{froese2013SINUM, froese2011JCP, oberman2015JCP} suitably combine a monotone scheme with a high-order scheme, in order to get a high-order approximation minimizing spurious oscillations. In this case the coupling is restricted to Eulerian schemes.

\item The Particle-in-Cell scheme \cite{birdsallPIC}, often used in plasma physics, combines a Lagrangian scheme with an Eulerian grid. Roughly speaking, Lagrangian particles are independently tracked, but an underlying grid is used to recover averaged quantities. At each time step, those quantities are used for additional computations which finally give the force acting on each particle. Then the procedure is repeated for the next time step.

\item The Particle Level-Set (PLS) method \cite{PLS} is a modification of the standard level-set method \cite{LS} for interfaces tracking, introduced to efficiently reduce the numerical diffusion appearing at regions of the evolving front with high curvature.  Lagrangian markers are posed in a neighborhood of the interested region, and moved independently by means of a high-order ODE solver, according to the characteristic flow of the level-set equation. 
Whenever a particle escapes from the considered neighbourhood, it means that level-set function describing the front starts to be smeared by the diffusion. Escaped particles are then employed to recover the mass loss, i.e., to locally reinitialize the level-set function. 

\item The Smoothed-Particle-Hydrodynamics (SPH) method \cite{SPH1, SPH2} (see also the Godunov-SPH \cite{puri2014JCP, inutsuka2002JCP}), is an efficient computational method proposed to give some regularity to a purely Lagrangian description of a phenomenon. This is done introducing suitable convolution kernels with compact supports and thinking each Lagrangian marker as a time-dependent Dirac measure. It follows that physical macroscopic quantities of interest can be obtained by simply summing the contribution of the sole particles within the cut-off radius of the kernel. Moreover, the spatial derivative of a quantity is easily computed by linearity using the gradient of the kernel. 

\end{itemize}

Finally, let us also mention the classical (W)ENO schemes \cite{carlini2005SISC, harten1987JCP, liu1994JCP, shu1997NASA} (see also \cite[Ch.\ 3]{falcone2014SIAM}), which provide high-order accuracy of the numerical solution and can be profitably used as main ingredients of the proposed scheme.

\subsubsection*{Paper organization} 
In Section \ref{sec:blendedscheme} we introduce the blended scheme for equations \eqref{AE} and \eqref{CL}, while in 
Section \ref{sec:theory} we study its theoretical properties. 
In Section \ref{sec:modifiedequation} we compute the modified equation and we focus on some  particular cases. 
In Section \ref{sec:num_tests} we present several numerical tests in order to enlighten the features of our approach. 
Finally, we sketch some conclusions and perspectives.

%%%%%%%%%%%%%%%%%%%%%%%%%%%%%%%%%%%%%%%
%%%%%%%%%%%%%%%%%%%%%%%%%%%%%%%%%%%%%%%
%%%%%%%%%%%%%%%%%%%%%%%%%%%%%%%%%%%%%%%
\section{Blended numerical scheme}\label{sec:blendedscheme}
To avoid cumbersome notations, let us assume $d=1$, results being easily generalized to any dimension. In addition, let us consider the coupling of only two numerical schemes. 
Let $\Omega^T:=\Omega\times[0,T]$ be the domain in which the solution will be computed, with $\Omega\subset\R$ a bounded domain and $T>0$ the final time. Let us introduce a structured grid $\GG$ in $\Omega^T$ and denote by $\{(x_i,t^n)\}$, $i=0,\ldots,N_C-1,\ n=0,\ldots,N_T-1$ the grid nodes, by $N_C$, $N_T$, respectively, the number of space and time nodes, and by $\Dx$, $\Dt$, respectively, the space and time discretization steps. 
The space grid cell is $C_i=\big[x_i-\frac{\Dx}{2},x_i+\frac{\Dx}{2}\big)$.

\subsection{General idea}
Let us define the projection on the grid of the exact solution $u$ of \eqref{AE}
\begin{equation}\label{defU}
U_i^n:=\frac{1}{\Dx}\int_{C_i}u(x,t^n)dx,\qquad i=0,\ldots,N_C-1,\qquad n=0,\ldots,N_T-1,
\end{equation}
and in particular the discrete initial condition
\begin{equation}\label{IC}
U_i^0=\frac{1}{\Dx}\int_{C_i}\bar u(x)dx,\qquad i=0,\ldots,N_C-1.
\end{equation}
Let us also denote respectively by $\{\W_i^n\}_{i,n}$, $\{\V_i^n\}_{i,n}$ two discrete solutions to \eqref{AE} at $(x_i,t^n)$, obtained by two (given) explicit-in-time numerical schemes of the form
\begin{equation}\label{PDEsolvers}
\left\{
\begin{array}{ll}
\W^{n+1}=\S_1[\W^{n}]\\ [2mm]
\W^0=U^0
\end{array}
\right.,
\qquad
\left\{
\begin{array}{ll}
\V^{n+1}=\S_2[\V^{n}]\\ [2mm]
\V^0=U^0
\end{array}
\right.,
\qquad
n=0,\ldots,N_T-1,
\end{equation}
where $\S_1,\S_2:\R^{N_C}\to\R^{N_C}$, 
$\W^n:=(\W_0^n,\ldots,\W_{N_C-1}^n)$ and $\V^n:=(\V_0^n,\ldots,\V_{N_C-1}^n)$. The treatment of boundary conditions is included in the definitions of the $\S_i$'s.

\medskip

We introduce now the \emph{blended scheme} defining the two new discrete solutions $\{W_i^n\}_{i,n}$ and $\{V_i^n\}_{i,n}$ to \eqref{AE}, with $W^n:=(W_0^n,\ldots,W_{N_C-1}^n)$ and $V^n:=(V_0^n,\ldots,V_{N_C-1}^n)$, by coupling $\S_1$ and $\S_2$ as follows 
\begin{equation}\tag{BS}\label{BlendedScheme}
\boxed{
\left\{
\begin{array}{l}
W^{n+1}=\thetauno \S_1[W^n] + (1-\thetauno) \S_2[V^n], \\ [3mm]
V^{n+1}=(1-\thetadue) \S_1[W^n] + \thetadue \S_2[V^n],
\end{array}
\right .
\qquad n=0,\ldots,N_T-1
}
\end{equation}
for some $\thetauno,\thetadue\in[0,1]$ and initial conditions $W^0=V^0=U^0$ as defined in \eqref{IC}. The \textit{coupling parameters} $\thetauno,\thetadue$ manage the convex combination between the schemes. Some remarks are in order:
\begin{enumerate}
\item[(i)] if $(\thetauno,\thetadue)\neq(1,1)$ the coupling occurs \emph{at each time step}, hence we do not get just a trivial interpolation between the two independent numerical approximations. On the contrary, it is quite hard to say which behavior will be actually observed in the end. Moreover, there arises a natural question about which couple $(\thetauno^*,\thetadue^*)$ would lead to the best result. We will try to give a partial answer to this question; 
\item[(ii)] if $(\thetauno,\thetadue)=(1,1)$, the two schemes evolve independently and $W^n\equiv \W^n$, $V^n\equiv \V^n$; 
\item[(iii)] if either $\thetauno=1$ or $\thetadue=1$, the blend is \emph{one-way} and only one solution is affected by the other; 
\item[(iv)] if $(\thetauno,\thetadue)=(0,0)$ the new scheme alternates the two original schemes at every time step, starting from $\S_2$ in the case of $W^n$ and $\S_1$ in the case of $V^n$;
\item[(v)] if $\thetauno\in(0,1)$ and $\thetadue=0$, \eqref{BlendedScheme} reduces to a single two-step scheme
\begin{equation*}
\begin{array}{ll}
W^1=\thetauno\S_1[U^0]+(1-\thetauno)\S_2[U^0] & \\ [3mm]
W^{n+1}=\thetauno\S_1[W^n]+(1-\thetauno)\S_2[\S_1[W^{n-1}]], & n\geq 1;
\end{array}
\end{equation*}
%This case is particularly important because numerical evidence tells us that the optimal blend is sometimes of this kind;
%
\item[(vi)] if $\thetauno=1-\thetadue$ we have $W^n\equiv V^n$ and \eqref{BlendedScheme} reduces to a single scheme
\begin{equation*}
W^{n+1}=\thetauno \S_1[W^n] + (1-\thetauno) \S_2[W^n], \qquad n\geq 0;
\end{equation*}
\item[(vii)] generalization to multi-blend is straightforward. If we have $\mathcal{L}$ schemes $\S_1,\ldots,{\S}_{\mathcal{L}}$, the blended scheme is 
$$
\left(
\begin{array}{c}
W^{n+1,(1)}\\
\vdots\\
W^{n+1,(\mathcal{L})}
\end{array}
\right)
=
\left(
\begin{array}{ccc}
\lambda_1^1 & \ldots & \lambda_\mathcal{L}^1\\
\vdots & \ddots & \vdots\\
\lambda_1^\mathcal{L} & \ldots & \lambda_\mathcal{L}^\mathcal{L}
\end{array}
\right)
\left(
\begin{array}{c}
\S_1[W^{n,(1)}]\\
\vdots\\
\S_\mathcal{L}[W^{n,(\mathcal{L})}]
\end{array}
\right)
$$
where the coupling parameters $\{\lambda_i^j\}$ are such that $\sum_{j=1}^\mathcal{L} \lambda_j^i=1$, for all $i=1,\ldots,\mathcal{L}$;
\item[(viii)] a more sophisticated blend could allow $(\thetauno,\thetadue)$ to depend on space, time, and the solution itself. %Note that in the space-dependent case mass conservation is no longer assured. 
In this way we can construct \emph{ad hoc} recipes wherever the solution shows some critical behaviour.
\end{enumerate}
\subsection{Multiscale coupling}\label{sec:multiscalecoupling}
In this section we transform the blended scheme \eqref{BlendedScheme} in a multiscale scheme, by choosing an Eulerian approximation for $\S_1$ (e.g., Upwind, Lax-Friedrichs, (W)ENO, Godunov, etc.) and a Lagrangian approximation for $\S_2$.
First, we consider a number of Lagrangian particles which move in $\Omega$ driven by the velocity field $A$. We assume that, at the initial time $t=0$, we have $\Np$ particles. We denote the (positions of the) particles at time $t$ by $\Pp(t)$, for $\alpha=0,\ldots,\Np-1$.  We also denote by $\Nppc=\frac{\Np}{N_C}\in\mathbb Q$ the number of particles per cell.
\subsubsection*{Initialization}
At time $t=0$, the particles are uniformly placed in $\Omega$, i.e.,
\begin{equation}\label{first_sampling}
%\P_\alpha(0)=\Pip(0)=x_i-\frac{\Dx}{2}+\left( p+\frac12 \right) \frac{\Dx}{\Nppc} %+p \frac{\Dx}{\Nppc}
%, \qquad i=0,\ldots,N_C-1, \quad p=0,\ldots,\Nppc-1,
\Pp(0)=x_0+\alpha\Delta p,\quad \alpha=0,\ldots,\Np-1\,, \qquad \Delta p=\frac{x_{N_C-1}-x_0}{\Np-1}.
\end{equation}
Moreover, each particle is endowed with a ``weight'', i.e.\ a portion of the mass transported by \eqref{AE}. Denoting by $\M_\alpha(t)$ the mass assigned to the $\alpha$-th particle at time $t$, we set
\begin{equation}\label{rectangleODE}
\M_\alpha(0)=\bar u\big(\P_\alpha(0)\big)\frac{\Dx}{\Nppc}.
\end{equation}
Note that, in a purely Lagrangian framework, there is no need to consider a time-dependent mass for each particle because it remains constant during the evolution. Therefore, one can initially distribute the particles in space according to $\bar u$, interpreting the initial condition as the probability density function of the mass. Here instead we define the mass $\M_\alpha$ as a time-dependent variable because of the combination with the Eulerian part of the algorithm, which will modify the masses at each time step. 

Note that one can avoid to put particles in regions in which $\bar u=0$, thus saving CPU time.  For simplicity, we put particles all along $\Omega$, regardless of the values of $\bar u$. 
\subsubsection*{Particles evolution}
As already mentioned in the introduction, the dynamics of each particle $\alpha$ is given by the ODE
\begin{equation}\label{ODEs}
\dot{\P}_\alpha(t)=A(\P_\alpha(t)).
\end{equation}
It is then natural to solve \eqref{ODEs} by means of an ODE solver. Focusing on one-step ODE solvers, the general algorithm has the form
\begin{equation}\label{ODEsolver}
\left\{
\begin{array}{l}
\P_\alpha^{n+1}=\ODEsolver(\P_\alpha^n), \qquad n=0,\ldots,N_T-1, \\ [2mm]
\P_\alpha^0=\P_\alpha(0),
\end{array}
\right.
\end{equation}
where $\P_\alpha^n:=\P_\alpha(t^n)$.
\subsubsection*{Computation of the Lagrangian density}
In order to combine the approximate Eulerian density with the Lagrangian computation we need a suitable averaging of the masses of the particles $\P_\alpha$'s, so to recover an approximate Lagrangian density. To this end we set 
\begin{equation}\label{S2micro}
\widehat V^{n+1}_i:=(\S_2[\P^{n+1},\M^n])_i=
\frac{1}{\Dx}\sum_{\alpha\ \! :\ \! \P_\alpha^{n+1}\in C_i} \M_\alpha^n.
\end{equation}
where $\M_\alpha^n=\M_\alpha(t^n)$. 
Let us explain this definition: at any time $t^n$ and for any cell $C_i$, we find the particles falling in the cell $C_i$, summing all their masses and dividing by the area of the cell.
Note that $\widehat V^{n+1}$ is computed by the updated positions of the particles $\P^{n+1}$ and the last available values for the masses.

\subsubsection*{Blending and mass update}
Given $\widehat W^{n+1}_i:=(\S_1[W^n])_i$ and $\widehat V^{n+1}_i$ as in \eqref{S2micro} at any space and time node, we are ready to run the scheme \eqref{BlendedScheme} and compute the blended approximations $W$ and $V$. However, modifying the Lagrangian density (from $\widehat V$ to $V$) leaves the particles  unchanged (both masses $\M$ and positions $\P$). Since the Lagrangian density is computed \emph{at each time step} by means of $\M$ and $\P$, see \eqref{S2micro}, we are going to lose the coupling effect after every time step. To fix this, we adopt the following strategy: at every time step, we modify the masses $\M_\alpha^n$ according to the modified Lagrangian density $V^{n+1}$. More, precisely, we denote by $\Lambda_i^n$ the number of particles falling in the cell $C_i$ at time $t^n$ and update their mass by setting:
\begin{equation}\label{massupdate}
\M^{n+1}_\alpha=\M^{n}_\alpha+
\left\{
\begin{array}{ll}
\displaystyle\frac{\Dx}{\Lambda_i^n}\Big(V_i^{n+1}-(\S_2[\P^{n+1},\M^n])_i\Big), & \text{if } \Lambda_i^n>0 \\ [4mm]
0, & \text{otherwise.}
\end{array}
\right .
\end{equation}
Obviously, since the Lagrangian density is constant on each cell, all the particles falling in the same cell are equally affected. Note that, as opposed to the masses, the positions $\P_\alpha$ are not affected by the coupling.

\subsubsection*{Final algorithm}
Summarizing, the complete algorithm for one time step $(W^n,V^n;\P^n,\M^n)$ $\to$ $(W^{n+1},V^{n+1};\P^{n+1},\M^{n+1})$ is the following:
\begin{enumerate}
\item Compute $\widehat W^{n+1}_i=(\S_1[W^n])_i$.
\item Compute $\P^{n+1}_\alpha$ by \eqref{ODEsolver} (all $\alpha$'s).
\item Compute $\widehat V^{n+1}_i=(\S_2[\P^{n+1},\M^n])_i$ by \eqref{S2micro} (all $i$'s).
\item Compute $W^{n+1}$ and $V^{n+1}$ by \eqref{BlendedScheme}.
\item Compute $\M^{n+1}_\alpha$ by \eqref{massupdate} (all $\alpha$'s).
\end{enumerate}

\begin{rmrk} (Partial coupling)
If $\thetadue=1$ the algorithm simplifies since there is no need to correct the masses, see \eqref{massupdate}.
\end{rmrk}
\begin{rmrk} (Parallelization)
The Lagrangian part of the multiscale algorithm can be easily and efficiently parallelized on both distributed and shared memory architectures since particles do not need to communicate with each other. In particular it is not needed to find the neighboring particles of each particle (as required in, e.g., SPH method). 
\end{rmrk}
\subsection{Extension to one-dimensional nonlinear conservation laws}\label{sec:fullynonlinear}
In this section we show how the blended scheme can be extended to equation \eqref{CL}. The Eulerian-Eulerian coupling does not need any modification and then it is no further discussed. The multiscale coupling (Section \ref{sec:multiscalecoupling}) needs instead some corrections. 
In conservation laws, the velocity $A$ of the particles (not that of the characteristic curves) can be recovered from the flux $f$ as $A(u):=\frac{f(u)}{u}$ (if $u\neq 0$) and equation \eqref{ODEs} must be replaced by
$$
\dot\P_\alpha(t)=A\big(u(\P_\alpha(t),t)\big).
$$
This means that the Lagrangian scheme needs the density $u$ to update the positions of the particles. The blended scheme provides at each time step two approximate densities, namely $W^n$ and $V^n$. Hence we solve either
\begin{equation}\label{fnl_WoV}
\dot\P_\alpha(t^n)=A(W^n_{j_{\alpha,n}}) \qquad \text{or} \qquad \dot\P_\alpha(t^n)=A(V^n_{j_{\alpha,n}})\,,
\end{equation} 
where $j_{\alpha,n}$ is the cell where the $\alpha$-th particle falls in at time $t^n$. Note that, choosing $W^n$, we introduces an additional coupling: the Lagrangian scheme now depends on the Eulerian scheme even for an uncoupled evolution $(\thetauno,\thetadue)=(1,1)$.  
\subsection{Estimation of coupling parameters $(\thetauno,\thetadue)$ by Richardson extrapolation}\label{sec:richardson}
The blended scheme \eqref{BlendedScheme} is actually a \textit{family} of numerical schemes indexed by two coupling parameters $(\thetauno,\thetadue)\in[0,1]^2$. Generally speaking, we expect that only some couples lead to a solution $W$ or $V$ which is better than both solutions $\W$ and $\V$ obtained by the uncoupled schemes. The question arises how to find an advantageous couple \textit{a priori}, before running the computation and without the knowledge of the exact solution. In some simple cases one can compute the optimal parameters analytically, but in general this is impossible. In order to provide a general method, we adopt the following strategy.

First of all, let us denote the solutions of the blended scheme \eqref{BlendedScheme} by $\{W^n_i(\thetauno,\thetadue)\}_{i,n}$ and $\{V^n_i(\thetauno,\thetadue)\}_{i,n}$ in order to make it explicit the dependence on the coupling parameters. Moreover, let us define the $L^1$ errors at final time as
\begin{equation}\label{def:E1}
\begin{split}
E^1[W(\thetauno,\thetadue)]:=\sum_i \Big|W^{N_T}_i(\thetauno,\thetadue)-U^{N_T}_i\Big|\Dx,\qquad
E^1[V(\thetauno,\thetadue)]:=\sum_i \Big|V^{N_T}_i(\thetauno,\thetadue)-U^{N_T}_i\Big|\Dx,
\end{split}
\end{equation}
and define the two optimal $L^1$ couplings as
\begin{equation}\label{def:theta*}
\begin{split}
(\thetauno^*,\thetadue^*)_W:=\arg\min_{\thetauno,\thetadue} E^1[W(\thetauno,\thetadue)],\qquad
(\thetauno^*,\thetadue^*)_V:=\arg\min_{\thetauno,\thetadue} E^1[V(\thetauno,\thetadue)].
\end{split}
\end{equation}
Following the classical idea of Richardson extrapolation, we can run the blended scheme on two (coarse) grids, $\GG'$ and $\GG''$, thus obtaining, for each $(\thetauno,\thetadue)$, four solutions $W(\thetauno,\thetadue;\GG')$, $V(\thetauno,\thetadue;\GG')$, $W(\thetauno,\thetadue;\GG'')$, and $V(\thetauno,\thetadue;\GG'')$. Then we define the two Richardson $L^1$ \textit{error indicators} as
\begin{equation}\label{def:deltaR}
\begin{split}
\delta^R_W(\thetauno,\thetadue;\GG',\GG''):=\sum_{i\in\GG'}\left|W_i^{N_T}(\thetauno,\thetadue;\GG')-W_{h(i)}^{N_T}(\thetauno,\thetadue;\GG'')\right|\Dx,\\
\delta^R_V(\thetauno,\thetadue;\GG',\GG''):=\sum_{i\in\GG'}\left|V_i^{N_T}(\thetauno,\thetadue;\GG')-V_{h(i)}^{N_T}(\thetauno,\thetadue;\GG'')\right|\Dx,
\end{split}
\end{equation}
where $h(i)$ gives the best correspondence from node $i$ in $\GG'$ to a node in $\GG''$. 
The estimates for $(\thetauno^*,\thetadue^*)$ are then given by
\begin{equation}\label{def:thetaR}
\begin{split}
(\thetauno^R,\thetadue^R)_{W,\GG',\GG''}:=\arg\min_{\thetauno,\thetadue} \delta^R_W(\thetauno,\thetadue;\GG',\GG''), \\
(\thetauno^R,\thetadue^R)_{V,\GG',\GG''}:=\arg\min_{\thetauno,\thetadue} \delta^R_V(\thetauno,\thetadue;\GG',\GG'').
\end{split}
\end{equation}
\begin{rmrk}
It is possible to avoid an exhaustive search in $[0,1]^2$ in \eqref{def:thetaR} by employing a descend method in the space $(\thetauno,\thetadue)$ starting from some initial guess $(\thetauno_0,\thetadue_0)$. One can also progressively refine the grid (in the space $(\thetauno,\thetadue)$) around the argmin found in the previous level.
\end{rmrk}
\begin{rmrk}
A necessary condition for the method to work is that $(\thetauno^*,\thetadue^*)$ is stable with respect to the grid size. In Section \ref{sec:num_tests} we will investigate this property. Moreover, we expect the method to be sufficiently fast if performed on very coarse grids, in particular much coarser than the final grid where the actual computation must be performed. 
\end{rmrk}

\section{Theoretical analysis}\label{sec:theory}

We address in this section the main theoretical issues for the blended scheme solving equation \eqref{AE} with $d=1$. Note that, rather than perform a complete convergence analysis (which apparently requires more technical subtleties), we will analyze the scheme in the two situations which seem to be more effective and that are actually used in numerical examples, i.e., the coupling of two Eulerian schemes and the use of a Lagrangian scheme to correct an Eulerian scheme in a multi-scale approach (i.e.\ $\thetauno\in[0,1)$ and $\thetadue=1$).

\subsection{Mass conservation}

We first prove that the scheme \eqref{BlendedScheme} is conservative, provided the two constitutive schemes are so. 
Let us use again the notation 
\begin{equation}\label{WeVcappuccio}
\widehat W_i^{n+1} = (\S_1[W^n])_i, \qquad\qquad
\widehat V_i^{n+1} = (\S_2[V^n])_i
\end{equation}
for a single time step of the uncoupled schemes (with a slight abuse of notation in the case of multiscale algorithm). 

\begin{prpstn}
Assume that the solvers \eqref{WeVcappuccio} are \emph{conservative}, i.e. 
\begin{equation}\label{teo:LAEconservativo_ipotesiPDEsolver}
\sum_i \widehat W_i^{n+1}=\sum_i W_i^n,
\qquad
\sum_i \widehat V_i^{n+1}=\sum_i V_i^n,
\end{equation}
and that 
\begin{equation}\label{ICuguali}
\sum_i W^0_i = \sum_i V^0_i=\sum_i U^0_i. 
\end{equation} 
Then, the scheme \eqref{BlendedScheme} is also conservative. 
\end{prpstn}
\begin{proof}
By definition of the scheme \eqref{BlendedScheme} and \eqref{WeVcappuccio}, we have 
\begin{equation}\label{teo:LAEconservativo_1}
\begin{cases}
\sum_i W_i^{n+1} = \thetauno \sum_i \widehat W_i^{n+1} + (1-\thetauno) \sum_i \widehat V_i^{n+1} \\
\sum_i V_i^{n+1} = (1-\thetadue) \sum_i \widehat W_i^{n+1} + \thetadue \sum_i \widehat V_i^{n+1},
\end{cases}
\end{equation}
so that, using the discrete conservation of both schemes \eqref{teo:LAEconservativo_ipotesiPDEsolver}, we have
$$
\begin{cases}
\sum_i W_i^{n+1} = \thetauno \sum_i W_i^n + (1-\thetauno) \sum_i V_i^n \\
\sum_i V_i^{n+1} = (1-\thetadue) \sum_i W_i^n + \thetadue \sum_i V_i^n,
\end{cases}
$$
which implies that
$$
\sum_i W_i^{n+1}, \sum_i V_i^{n+1} \in \left[\min\left(\sum_i W_i^n, \sum_i V_i^n\right), \ \max\left(\sum_i W_i^n, \sum_i V_i^n\right)\right].
$$
The proof is concluded by induction noting that \eqref{ICuguali} holds.
\end{proof}

\begin{rmrk}
The previous result remains true if both parameters $\thetauno$ and $\thetadue$ depend on the step index $n$, but not if they depend on the cell index $i$.
\end{rmrk}

\subsection{Coupling between two Eulerian schemes}

Hereafter we will assume for simplicity that the schemes $\S_1,\S_2$ are linear, so that we can write 
$$
\W^{n+1}=S_1 \W^n \quad,\qquad \V^{n+1}=S_2 \V^n,
$$
where $S_1, S_2$ are two $N_C\times N_C$ matrices.
Let us also assume that each of the schemes is stable, and in particular that they satisfy 
\begin{equation}\label{stab}
\| S_1 \|_{p_1} \le 1 + C_1 \Delta t \quad , \qquad
\| S_2 \|_{p_2} \le 1 + C_2 \Delta t,
\end{equation}
with constants $C_1, C_2$ independent of $\Delta x$ and $\Delta t$, and for some (possibly different) suitable norms $\|\cdot\|_{p_i}$ $(i=1,2)$ which could be required for the convergence analysis of the two schemes (e.g., $\infty$-norm for a monotone scheme and $2$-norm for a high order scheme). We also assume to use scaled $l^p$ norms defined by
\begin{equation}\label{norma}
\|U\|_p = \Dx^{1/p}\left(\sum_j |u_j|^p\right)^{1/p}.
\end{equation}
Now, the scheme \eqref{BlendedScheme} can be rewritten in matrix form as
\begin{equation}\label{LAE_matrix}
\begin{pmatrix} W^{n+1} \\ V^{n+1} \end{pmatrix}
%& = &
\begin{pmatrix} \thetauno I & (1-\thetauno) I \\ (1-\thetadue) I & \thetadue I \end{pmatrix}
\begin{pmatrix} S_1 & 0 \\ 0 & S_2 \end{pmatrix}
\begin{pmatrix} W^n \\ V^n \end{pmatrix} = %\nonumber 
%& = & 
\Theta S \begin{pmatrix} W^n \\ V^n \end{pmatrix}
\end{equation}
where $I$ denotes the $N_C\times N_C$ identity matrix, and
$$
\Theta =
\begin{pmatrix} \thetauno I & (1-\thetauno) I \\ (1-\thetadue) I & \thetadue I \end{pmatrix}
\quad, \qquad
S =\begin{pmatrix} S_1 & 0 \\ 0 & S_2 \end{pmatrix}.
$$

In order to analyse convergence for the scheme \eqref{BlendedScheme}, we define a suitable norm for the numerical solution $(W^n, V^n)$ in the product space $\R^{N_C}\times \R^{N_C}$ as
$$
\| (W, V) \|_B := \max(\|W\|_{p_1},\|V\|_{p_2}).
$$
Note that
\begin{equation}\label{Theta_norm}
\|\Theta\|_B = 1 \quad,\qquad \|S\|_B \le 1+\max\{C_1,C_2\}\Delta t,
\end{equation}
which implies that \eqref{BlendedScheme} is stable in the norm $\|\cdot\|_B$, since
$$
\|\Theta S\|_B \le \|\Theta\|_B \|S\|_B \le 1+\max\{C_1,C_2\}\Delta t.
$$
It is also immediate to verify that \eqref{BlendedScheme} is consistent, and therefore convergent. In fact, assume to work in the norm $\|\cdot\|_B$ and apply the usual definition of consistency for a discrete solution at time $t^n$ defined by
$$
\begin{pmatrix}
W^n \\ V^n
\end{pmatrix}
=
\begin{pmatrix}
U^n \\ U^n
\end{pmatrix},
$$
where $U^n$ is defined in \eqref{defU} for a smooth $u(x,t)$. Taken into account that the discretization parameters are $\Dt$, $\Dx$, the scheme is consistent if and only if
\begin{equation}\label{cons1}
\begin{pmatrix} U^{n+1} \\ U^{n+1} \end{pmatrix} = \Theta S \begin{pmatrix} U^n \\ U^n \end{pmatrix} + 
\Delta t \begin{pmatrix} \sigma_1(\Dx,\Dt) \\ \sigma_2(\Dx,\Dt) \end{pmatrix}
\end{equation}
with $\sigma_1,\sigma_2\to 0$ for $\Dt,\Dx\to 0$. Now, if the two schemes have respectively truncation errors $\tau_1$ and $\tau_2$, it is immediate to compute $\sigma_1$ and $\sigma_2$ as
$$
\begin{pmatrix} \sigma_1(\Dx,\Dt) \\ \sigma_2(\Dx,\Dt) \end{pmatrix} = \Theta \begin{pmatrix} \tau_1(\Dx,\Dt) \\ \tau_2(\Dx,\Dt) \end{pmatrix}
$$
which results in consistency for the scheme \eqref{BlendedScheme} if the two elementary blocks are consistent. We can therefore conclude with the following

\begin{thrm}
Let $S_1$, $S_2$ satisfy \eqref{stab}, and, for $i=1,2$,
$$
\tau_i(\Dx,\Dt) \to 0 \quad (\Dx,\Dt\to 0).
$$
Then, the blended scheme \eqref{BlendedScheme} is convergent in the norm $\|\cdot\|_B$ for any $\lambda,\mu\in[0,1]$.
\end{thrm}

\subsection{Convergence for the Lagrangian correction}

We recall that in this case we are assuming that the Lagrangian part of the scheme provides a correction for the Eulerian part. This requires to set $\mu=1$, whereas $\lambda\in [0,1]$, and in this case we clearly have $V^n = \V^n$. In order to have a more explicit notation, we will denote here the two components of the error at the $n$-th time step for the Eulerian and the Lagrangian part as $\epsilon_E^n$ and $\epsilon_L^n$,  respectively.

We start by computing the error for the pure Lagrangian part of the scheme, then we turn to the coupled scheme.

\subsubsection*{Step 1. Bound on $\epsilon_L^n:=U^n-V^n$}
\begin{figure}[h!]
\begin{center}
\begin{tabular}{c}
\includegraphics[width=0.32\textwidth]{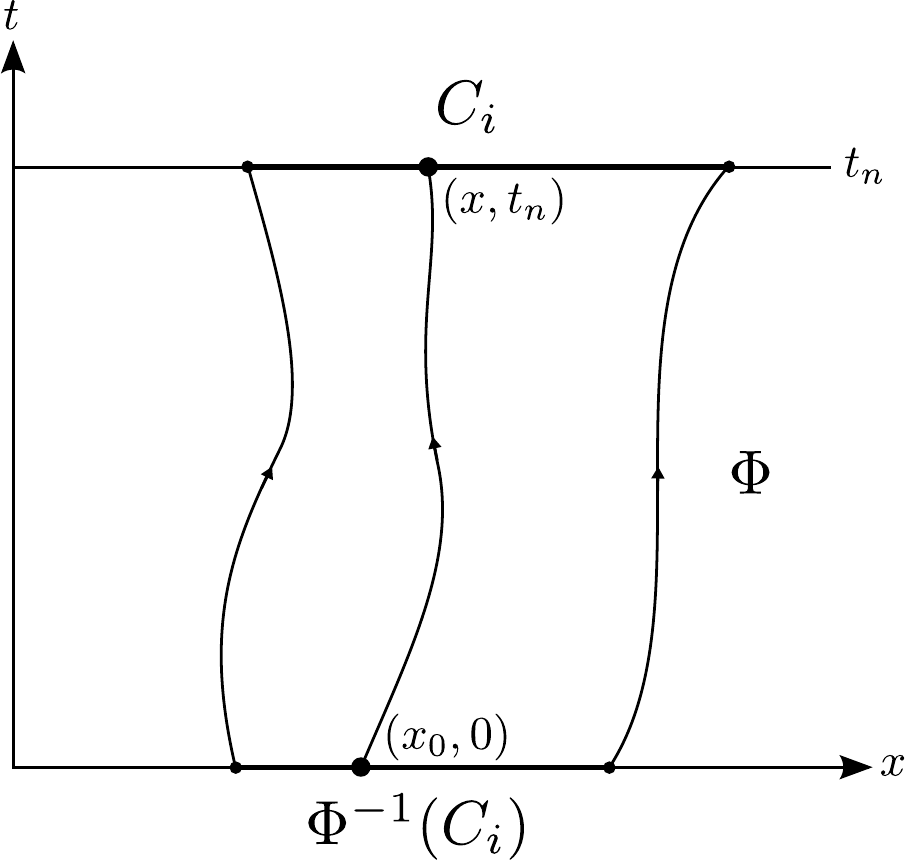} 
\end{tabular}
\end{center}
\caption{The characteristic map $\Phi$.}
\label{fig:Map}
\end{figure}

We denote by $\Phi: (x_0,0)\mapsto (x,t_n)$ the mapping between the foot $x_0$ of the \emph{exact} characteristic curve at time $t=0$ reaching the point $x$ at time $t_n$ (see Fig.\ \ref{fig:Map}). 

Similarly, we denote by $\Phi_\Delta$ the mapping corresponding to the \emph{discrete} characteristic curve, given by the ODE solver \eqref{ODEsolver}. 
We will make the standing assumption that this is a monotonic and invertible map (i.e., that no crossing between discrete characteristics occurs). In usual ODE solvers, this happens for $\Delta t$ sufficiently small. We also assume that the error of ODE scheme is of order $O(\Delta t^\gamma)$, for some $\gamma\geq 1$.

Then mass conservation and $u(x,0)=\bar u(x)$ imply
$$U^n_i=\frac{1}{\Delta x}\int_{C_i}u(x,t_n)dx=\frac{1}{\Delta x}\int_{\Phi^{-1}(C_i)}\bar u(y)dy=\frac{1}{\Delta x}\int_{C_i}\bar u(\Phi^{-1}(x))(\Phi^{-1}(x))^\prime dx\,,$$
where the last equality follows applying the change of variable $y=\Phi^{-1}(x)$.

We also define 
$$\tilde U^n_i=\frac{1}{\Delta x}\int_{C_i}\bar u(\Phi^{-1}_\Delta(x))(\Phi^{-1}_\Delta(x))^\prime dx$$
to split the estimate as $|U^n_i-V^n_i|\le |\tilde U^n_i- V^n_i|+|\tilde U^n_i - U^n_i|$.\\

Let us recall the definition of $V_i^n$ as in \eqref{S2micro},
$$V^n_i=\hat V^n_i=\frac{1}{\Dx}\sum_{\alpha\ \! :\ \! \P_\alpha^{n}\in C_i} \M_\alpha^n\,,$$
where the masses are computed by backtracking the Lagrangian particles up to the initial time.
Let $N$ be the number of particles $P_k^n$ falling in $C_i$, indexed for simplicity by $k=1,...,N$. Note that $N=O(N_{PC})$ by the Gronwall's lemma. Since $P_k^n$ is the evolution at time $t_n$ of the corresponding $P_k^0$ 
at time zero, we have
$$V^n_i=\frac{1}{\Dx}\sum_{k=1}^N \bar u(P_k^0)\frac{\Delta x}{N_{PC}}\,.$$
Using the discrete map $\Phi_\Delta$, we can decompose the cell $C_i$ as the union of the pre-images of uniform intervals centered at the particles, namely 
$$C_i=\bigcup_{k=1}^N\Phi_\Delta\left(\left[P_k^0-\frac{\Delta x}{2 N_{PC}},P_k^0+\frac{\Delta x}{2 N_{PC}}\right)\right)=:\bigcup_{k=1}^N\Phi_\Delta\left(\left[P_{k-\frac12}^0,P_{k+\frac12}^0\right)\right)\,.$$
Applying the change of variable $y=\Phi_\Delta^{-1}(x)$ we obtain
$$|\tilde U^n_i- V^n_i|\le \frac{1}{\Delta x}\sum_{k=1}^N \left|\int_{P_{k-\frac12}^0}^{P_{k+\frac12}^0} \bar u(y)dy- \bar u(P_k^0)\frac{\Delta x}{N_{PC}}\right|=
\frac{1}{\Delta x}\sum_{k=1}^N \frac{\Delta x}{N_{PC}}\left|\bar u\left(P_k^0+O\left(\frac{\Delta x}{N_{PC}}\right)\right)- \bar u(P_k^0)\right|\,,$$
where the last step follows from the mean value theorem and the fact that $P_{k+\frac12}^0 -P_{k-\frac12}^0=\frac{\Delta x}{N_{PC}}$.
Then we get
\begin{equation}\label{simoA}
|\tilde U^n_i- V^n_i|\le \frac{O(N_{PC})}{N_{PC}} L_{\bar u}O\left(\frac{\Delta x}{N_{PC}}\right)\le C\frac{\Delta x}{N_{PC}}\,,
\end{equation}
where $L_{\bar u}$ is the Lipschitz constant of the initial datum.

\medskip

\begin{figure}[h!]
\begin{center}
\begin{tabular}{c}
\includegraphics[width=0.32\textwidth]{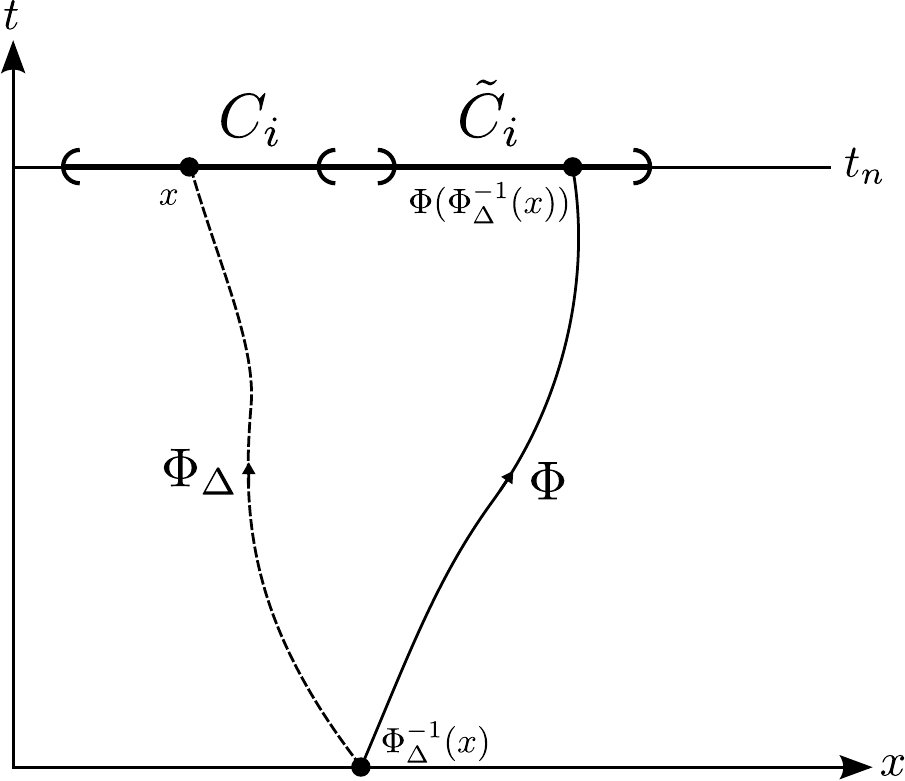} 
\end{tabular}
\end{center}
\caption{Combining exact and discrete characteristic maps $\Phi$ and $\Phi_\Delta$.}
\label{fig:MapErr}
\end{figure}

Now, by the change of variable $\Phi_\Delta^{-1}(x)=\Phi^{-1}(y)$ we get
$$ |\tilde U^n_i- U^n_i|\le \frac{1}{\Delta x} \left|\int_{\tilde C_i}\bar u(\Phi^{-1}(y))(\Phi^{-1}(y))^\prime dy-\int_{C_i}\bar u(\Phi^{-1}(x))(\Phi^{-1}(x))^\prime dx\right|\,,$$
where $\tilde C_i=\Phi(\Phi_\Delta^{-1}(C_i))$. In practice (see Fig.\ \ref{fig:MapErr}), we trace back each point $x\in C_i$ along the discrete characteristic
up to $\Phi_\Delta^{-1}(x)$, then we move forward along the exact characteristic up to $y=\Phi(\Phi_\Delta^{-1}(x))\in \tilde C_i$. 
\\
Since $x$ and $y$ share the same foot, the error estimate 
for the ODE solver provides a function $\phi$ such that
$$y=x+\phi(x)O(\Delta t^\gamma)\,,$$
namely a mapping between $\tilde C_i$ and $C_i$, which in turn gives 
$$|\tilde U^n_i- U^n_i|\le \frac{1}{\Delta x} \int_{C_i}\Big|\bar u(\Phi^{-1}(x+\phi(x)O(\Delta t^\gamma)))(\Phi^{-1}(x+\phi(x)O(\Delta t^\gamma)))^\prime -\bar u(\Phi^{-1}(x))(\Phi^{-1}(x))^\prime \Big|dx + O(\Delta t^\gamma)\le$$
\begin{equation}\label{simoB}
\le  C \frac{1}{\Delta x}\int_{C_i}|\phi(x)|O(\Delta t^\gamma) + O(\Delta t^\gamma)\le C \Delta t^\gamma\,.
\end{equation}

Combining \eqref{simoA} with \eqref{simoB}, for any norm in the form \eqref{norma}, the following bound holds:
$$\|\epsilon_L^n\|_p \le C\left(\Delta t^\gamma+\frac{\Delta x}{N_{PC}}\right)\,.$$

\subsubsection*{Step 2. Bound on $\epsilon_E^n:=U^n-W^n$}

Writing \eqref{BlendedScheme} between $n-1$ and $n$, and substituting the second relationship in the first one, we get
$$
W^n = \lambda S_1W^{n-1} + (1-\lambda)\left(U^n-\epsilon_L^n\right),
$$
which gives in turn, for the error on the component $W^n$,
\begin{eqnarray}\label{eps_1}
\epsilon_E^n = U^n-W^n & = & \lambda \big[U^n - S_1W^{n-1}\big] + (1-\lambda)\epsilon_L^n = \nonumber \\
& = & \lambda \big[\left(U^n-S_1U^{n-1}\right)+\left(S_1U^{n-1}-S_1W^{n-1}\right)\big] + (1-\lambda)\epsilon_L^n.
\end{eqnarray}
Using now the consistency error $\tau_1$ of the Eulerian part of the scheme, we can write the first term in square brackets as
\begin{equation}
U^n-S_1U^{n-1} = \Dt\, \tau_1(\Dx,\Dt)
\end{equation}
which gives, when used in \eqref{eps_1}:
\begin{equation}\label{stima2}
\epsilon_E^n = \lambda \big[\Dt\, \tau_1(\Dx,\Dt)+\left(S_1(U^{n-1}-W^{n-1})\right)\big] + (1-\lambda)\epsilon_L^n.
\end{equation}
Passing to the norms, we obtain therefore
\begin{equation}\label{stima3}
\left\|\epsilon_E^n\right\|_p \le \lambda (1+C_1\Dt) \left\|\epsilon_E^{n-1}\right\|_p + \lambda\Dt \left\|\tau_1(\Dx,\Dt)\right\|_p + (1-\lambda)\left\|\epsilon_L^n\right\|_p.
\end{equation}
Defining now $\zeta := \lambda (1+C_1\Dt)$ and iterating back the estimate, we obtain
\begin{eqnarray}\label{stima_eps1}
\left\|\epsilon_E^n\right\|_p & \le & \zeta \left\|\epsilon_E^{n-1}\right\|_p + \lambda\Dt \left\|\tau_1(\Dx,\Dt)\right\|_p + (1-\lambda)\left\|\epsilon_L^n\right\|_p \le \nonumber \\
& \le & \zeta^2 \left\|\epsilon_E^{n-2}\right\|_p + \lambda\Dt \left\|\tau_1(\Dx,\Dt)\right\|_p(1+\zeta) + (1-\lambda)\left(\left\|\epsilon_L^n\right\|_p + \zeta \left\|\epsilon_L^{n-1}\right\|_p\right) \le \cdots \le \nonumber \\
& \le & \zeta^n \left\|\epsilon_E^0\right\|_p + \lambda\Dt \left\|\tau_1(\Dx,\Dt)\right\|_p \sum_{k=0}^{n-1}\zeta^k + (1-\lambda)\sum_{k=0}^{n-1}\zeta^k\left\|\epsilon_L^{n-k}\right\|_p \le \nonumber \\
& \le & \zeta^n \left\|\epsilon_E^0\right\|_p + \left(\lambda\Dt \left\|\tau_1(\Dx,\Dt)\right\|_p + (1-\lambda)\max_k\left\|\epsilon_L^{n-k}\right\|_p \right) \sum_{k=0}^{n-1}\zeta^k.
\end{eqnarray}
Taking into account the definition of $\zeta$, the last term in \eqref{stima_eps1} can be bounded as
$$
\sum_{k=0}^{n-1}\zeta^k \le
\begin{cases}
\displaystyle \sum_{k=0}^\infty\zeta^k \le \frac{1}{1-(\lambda+\delta)} & \text{if } \lambda<1 \\ \\
\displaystyle \frac{e^{C_1T}-1}{C_1\Dt} & \text{if } \lambda=1,
\end{cases}
$$
where the first case clearly holds for $\Dt$ small enough to have $\lambda(1+C_1\Dt) \le \lambda+\delta < 1$. Dropping the trivial case $\lambda=1$ (which gives two uncoupled schemes in \eqref{BlendedScheme}), we obtain therefore, for $\lambda<1$ and $\Dt$ small enough, the estimate
$$
\left\|\epsilon_E^n\right\|_p \le \left\|\epsilon_E^0\right\|_p + C(\lambda) \left(\Dt \left\|\tau_1(\Dx,\Dt)\right\|_p + \max_k\left\|\epsilon_L^{n-k}\right\|_p \right).
$$

\subsubsection*{Step 3. Conclusions}

Finally, assuming that the initial condition has exact cell averages (i.e., that $\left\|\epsilon_E^0\right\|_p=0$), and recalling that $N_{PC}=\frac{\Np}{N_C}=\Np O(\Delta x)$, we obtain that the scheme converges and the error $\epsilon^n=(\epsilon_E^n,\epsilon_L^n)$ satisfies, in the norm $\|\cdot\|_B$, the error bound
\begin{equation}\label{StimaErr}
\left\|\epsilon^n\right\|_B \le 
C \left(\Dt \left\|\tau_1(\Dx,\Dt)\right\|_p + \Delta t^\gamma+\frac{\Delta x}{N_{PC}}\right)
=C \left(\Dt \left\|\tau_1(\Dx,\Dt)\right\|_p + \Delta t^\gamma+\frac{1}{\Np}\right),
\end{equation}
for any $n>0$ such that $t_n\in [0,T]$. We can then summarize this analysis in the following

\begin{thrm}
Let $S_1$ satisfy \eqref{stab}, and assume that $\tau_1(\Dx,\Dt) \to 0$ for $\Dx,\Dt\to 0$.

Then, for any $\lambda\in [0,1]$, the blended scheme \eqref{BlendedScheme}, as detailed in Subsection \ref{sec:multiscalecoupling}, is convergent in the norm $\|\cdot\|_B$ as $\Dx,\Dt\to 0,\ \Np\to+\infty$, and the estimate \eqref{StimaErr} holds. 
\end{thrm}

%%%%%%%%%%%%%%%%%%%%%%%%%%%%%%%%%%%%%%%
%%%%%%%%%%%%%%%%%%%%%%%%%%%%%%%%%%%%%%%
%%%%%%%%%%%%%%%%%%%%%%%%%%%%%%%%%%%%%%%
\section{Modified equations}\label{sec:modifiedequation}
In this section we derive a system of modified equations of \eqref{AE} (see \cite[Sect.\ 11.1]{levequebook}) associated to the scheme \eqref{BlendedScheme}, 
in order to give an interpretation of the effect of the coupling between the schemes from a continuous point of view. 

\subsection{The general case}
For our purposes, it is convenient to write first the scheme 
in the following implicit form, obtained by straightforward algebraic manipulations in the case $\thetauno\not = 1-\thetadue$: 
\begin{equation}\label{LAEbis}
\left\{
\begin{array}{l}
\displaystyle W^{n+1}-\S_1[W^n] =\frac{1-\thetauno}{\thetauno+\thetadue-1} (V^{n+1}-W^{n+1}) \\ [3mm]
\displaystyle V^{n+1}-\S_2[V^n] =\frac{1-\thetadue}{\thetauno+\thetadue-1} (W^{n+1}-V^{n+1})
\end{array}
\right .
\qquad n=0,\ldots,N_T-1\,.
\end{equation}
Now, let us assume $A(x)\equiv A$ and that, for a given $\Delta x/\Delta t$, 
the schemes $\mathcal{S}_1$ and $\mathcal{S}_2$ have respectively the modified equations
$$
w_t + Aw_x = \nu_1\Delta x^{p_1}\frac{\partial^{q_1}w}{\partial x^{q_1}}, \qquad 
v_t + Av_x = \nu_2\Delta x^{p_2}\frac{\partial^{q_2}v}{\partial x^{q_2}}, 
$$
for some $p_1, q_1, p_2, q_2$. Dividing by $\Delta t$ in \eqref{LAEbis}, 
using a Taylor expansion of the right hand sides up to first order in $\Delta t$ and 
restricting to the leading terms, we 
end up with the following system for $w$ and $v$:\\
\begin{equation}\label{modifeq_completa}
\left\{
\begin{array}{l}
\displaystyle w_t+Aw_x=\nu_1\Delta x^{p_1}\frac{\partial^{q_1}w}{\partial x^{q_1}}+\frac{1}{\Delta t}\frac{1-\thetauno}{\thetauno+\thetadue-1} (v-w)+\frac{1-\thetauno}{\thetauno+\thetadue-1}(v_t-w_t) \\ [3mm]
\displaystyle v_t+Av_x=\nu_2\Delta x^{p_2}\frac{\partial^{q_2}v}{\partial x^{q_2}}+\frac{1}{\Delta t}\frac{1-\thetadue}{\thetauno+\thetadue-1} (w-v)+\frac{1-\thetadue}{\thetauno+\thetadue-1}(w_t-v_t)
 \\\\w(x,0)=v(x,0)=\bar u(x)
\end{array}
\right .
\end{equation}

We clearly see that the dominant terms involving $1/\Delta t$ act as reaction terms, 
depending on the sign of the differences $\pm(v-w)$. 
%In the absence of the numerical dispersion (i.e. for $\nu_1=\nu_2=0$), since $w$ and $v$ coincide at the initial time, it easily follows that $w=v=u$ for all the times, where $u$ is the exact solution. 
%In the general case, even if we start with $w(x,0)=v(x,0)=\bar u(x)$, the time evolution of the solution to the system \eqref{modifeq_completa} is no longer trivial. 
We expect that, at the very beginning, $w$ and $v$ start to deviate from each other. This immediately activates both the reactions, driven by the weights $\frac{1-\thetauno}{\thetauno+\thetadue-1}$ and $\frac{1-\thetadue}{\thetauno+\thetadue-1}$ respectively, which 
push the solutions close to each other. 
%The final {\em product} is a pair of solutions which inherit features of both equations, according to the percentage of their mixing. 
%
\begin{rmrk}
The reaction terms appearing in the modified equation \ref{modifeq_completa} are similar to those appearing in \cite[Eq.(10)]{dimascio2016JCP}, where the authors combine the Lagrangian SPH solution with a Finite Volume scheme. There, a ``chimera''-like approach is used to force one solution towards the other, in order to minimize discontinuities and artifacts. 
\end{rmrk}

\subsection{The particular case $\thetauno=1-\thetadue$}\label{sec:particularcase3}

If $\thetauno=1-\thetadue$ we have $W^n\equiv V^n$ (cf.\ (vi) in Section \ref{sec:blendedscheme}), 
and with some minor manipulation we can rewrite (BS) as
$$
\frac{W^{n+1}-W^n}{\Delta t} = \lambda \frac{\mathcal{S}_1[W^n]-W^n}{\Delta t} + (1-\lambda) \frac{\mathcal{S}_2[W^n]-W^n}{\Delta t}.
$$
Restricting again to the leading terms, this yields the modified equation
$$
w_t + Aw_x = \left(\lambda\nu_1\Delta x^{p_1}\frac{\partial^{q_1}}{\partial x^{q_1}} + (1-\lambda)\nu_2\Delta x^{p_2}\frac{\partial^{q_2}}{\partial x^{q_2}} \right) w
$$
which shows a final dispersive term obtained by a mere convex combination of the two original ones.

Note that, if the two schemes are not of the same order, i.e., if $p_1\not=p_2$, then the scheme of lower order becomes asymptotically dominant if $\lambda$ and $\mu$ are kept constant along a refinement. On the other hand, if a nontrivial optimal combination exists, then both $\thetauno$ and $\thetadue$ are expected to depend on $\Delta x$.

\subsubsection*{Example 1: Lax-Wendroff + Beam-Warming} 
Assume that $\lambda=1-\mu$, and that we are coupling the Lax-Wendroff (LW) and the Beam-Warming (BW) schemes \cite[Chapt.\ 10]{levequebook} in the simple case $A(x)\equiv A>0$. As a consequence of having dispersive terms of the same order, but different sign,  these two schemes are known to have opposite behaviour with respect to numerical dispersion. In fact, the first one causes oscillations {\em behind} a discontinuity, while the second one {\em ahead} of it. We have therefore some chance to compensate the two opposite responses.

Denoting the Courant number by $\beta=A\Delta t/\Delta x$, we have for the LW scheme the modified equation
$$
w_t + Aw_x = -\frac{A\Delta x^2}{6}\left(1-\beta^2\right) w_{xxx}
$$
while the BW scheme gives
$$
w_t + Aw_x = -\frac{A\Delta x^2}{6}\left(2-3\beta+\beta^2\right) w_{xxx}.
$$
Now, the coupled scheme has a modified equation given by
$$
w_t + Aw_x = \frac{A\Delta x^2}{6}\left[\lambda\left(1-\beta^2\right)-(1-\lambda)\left(2-3\beta+\beta^2\right)\right] w_{xxx},
$$
and hence, imposing the term in square brackets to vanish, we obtain the condition 
$
3\lambda(1-\beta) = \beta^2-3\beta+2,
$
which finally gives the value of the optimal parameter
$
\lambda^* = \frac{2-\beta}{3}
$. 
Clearly, this procedure makes the {\em leading} dispersive term vanish, thus the modified equation of the blended scheme will now have a nonzero term in $w_{xxxx}$.
A straightforward computation shows that the resulting scheme is nothing but a third-order Upwind (or semi-Lagrangian) scheme obtained by interpolating the numerical solution $W^n$ at the foot of a characteristic line via a symmetric Lagrange polynomial of third degree.

\subsection{The particular case $\V\equiv U$}\label{sec:particularcase2}

Some additional considerations can be done in a even more special case, namely assuming that the scheme $\S_2$ corresponds to the {\em exact solution} (i.e.\ $\V\equiv U$):
$$\V_i^{n}=u(x_i,n\Delta t)=\bar u(x_i-An\Delta t)\,.$$
This analysis is interesting because gives some insights about the coupling with nondiffusive schemes. In particular, it shows the limit behavior of the coupling with a Lagrangian scheme when the number of particles tends to infinity.

Note that the blending of $\V$ with $\S_1$ does not make sense in this case, being $\V$ exact, therefore we keep $\thetadue$ fixed to $1$. Again, due to the numerical dispersion of $\S_1$, 
one could expect that the blended solution $W$ 
still shows a diffusion, just mitigated by the exact solution. 
Surprisingly, this is not what is observed in the numerical simulation, 
which instead shows that the numerical diffusion disappears after a short transient. This behavior can be explained by looking again at the modified equations. 
After a rearrangement of the terms containing $w_t$, the system \eqref{modifeq_completa} reduces to
\begin{equation}\label{eq_modif:upwind+esatta} 
\left\{
\begin{array}{l}
\displaystyle w_t+\thetauno A w_x=\lambda\nu\Delta x^{p}\frac{\partial^{q}w}{\partial x^{q}}+
\frac{1-\thetauno}{\Delta t}(u-w)+(1-\thetauno) u_t\\ [2mm]
u_t+Au_x=0\\ [2mm]
w(x,0)=u(x,0)=\bar u(x)
\end{array}
\right .
\end{equation}
for $\nu=\nu_1$, $p=p_1$ and $q=q_1$. 
It is easy to verify that the difference $r:=w-u$ satisfies 
\begin{equation}\label{eq:residuo}
r_t(x,t)+ \thetauno A r_x(x,t)= - \frac{1-\thetauno}{\Delta t} r(x,t) +\thetauno\nu \Delta x^{p}\frac{\partial^{q}}{\partial x^{q}}r(x,t)
+\thetauno\nu\Delta x^{p}\frac{\partial^{q}}{\partial x^{q}}u(x,t),
\end{equation}
with the initial condition $r(x,0)=0$.
Performing a Fourier transform with respect to the space variable, we get:
\begin{eqnarray}\label{eq:residuo2}
\dot{\hat r}(\xi,t) & = & - \thetauno A i\xi \hat r(\xi,t) - \frac{1-\thetauno}{\Delta t}\hat r(\xi,t) 
+\thetauno\nu \Delta x^{p} (i\xi)^{q} \hat r(\xi,t)+\thetauno\nu \Delta x^{p} (i\xi)^{q}\hat u(\xi,t) = \nonumber \\
& = & \left(- \thetauno A i\xi  -  \frac{1-\thetauno}{\Delta t} +\thetauno\nu \Delta x^{p} (i\xi)^{q}\right) 
\hat r(\xi,t)+ \thetauno\nu \Delta x^{p} (i\xi)^{q}\hat u(\xi,t),
\end{eqnarray}
where $\xi$ is the Fourier variable and $\hat r$ denotes the Fourier transform of $r$.
We observe now that \eqref{eq:residuo2} represents a family of decoupled linear ODEs indexed by $\xi$ 
with source terms given by $\thetauno\nu \Delta x^{p} (i\xi)^{q}\hat u(\xi,t)$, 
and that the coefficient in brackets has a negative real part for any $\xi$, provided the dispersion operator is dissipative (as it is usually the case in stable schemes).
Then, we infer that, for $t\to\infty$, the function $\hat r(\xi,t)$ (and hence, the function $r(x,t)$) 
remains bounded in $L^2$ provided $u(\cdot,t)$ has itself a bounded 2-norm. This shows that the numerical diffusion of the blended scheme is eventually balanced by the reaction term, and no further deteriorates the solution.
This observation is confirmed by the following numerical example. 
\subsubsection*{Example 2: Upwind + Exact}
Let us couple the Upwind scheme ($\S_1$) with an exact scheme ($\S_2$) with $\Omega=[0,20]$, $T=10$, $A=1$, $N_C=300$, $N_T=800$, and initial condition
\begin{equation}\label{baruSimone}
\bar u(x)=
\left\{
\begin{array}{ll}
\displaystyle\frac12 \left( 1+\cos(\pi(x-2))\right), & \text{ if } |x-2|<1, \\ [3mm]
0, & \text{ otherwise.}
\end{array}
\right.
\end{equation}
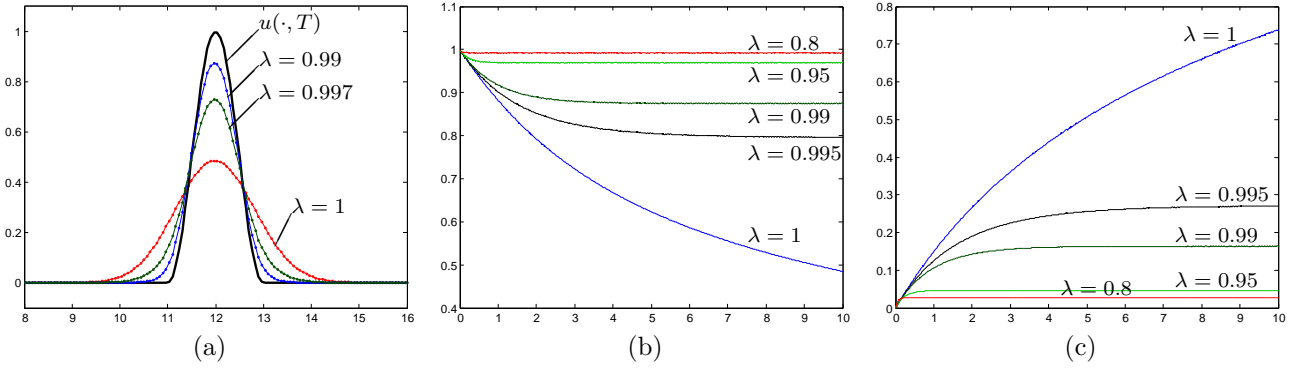
\begin{figure}[h!]
\begin{center}
\begin{tabular}{ccc}
\begin{overpic}
[width=0.32\textwidth]{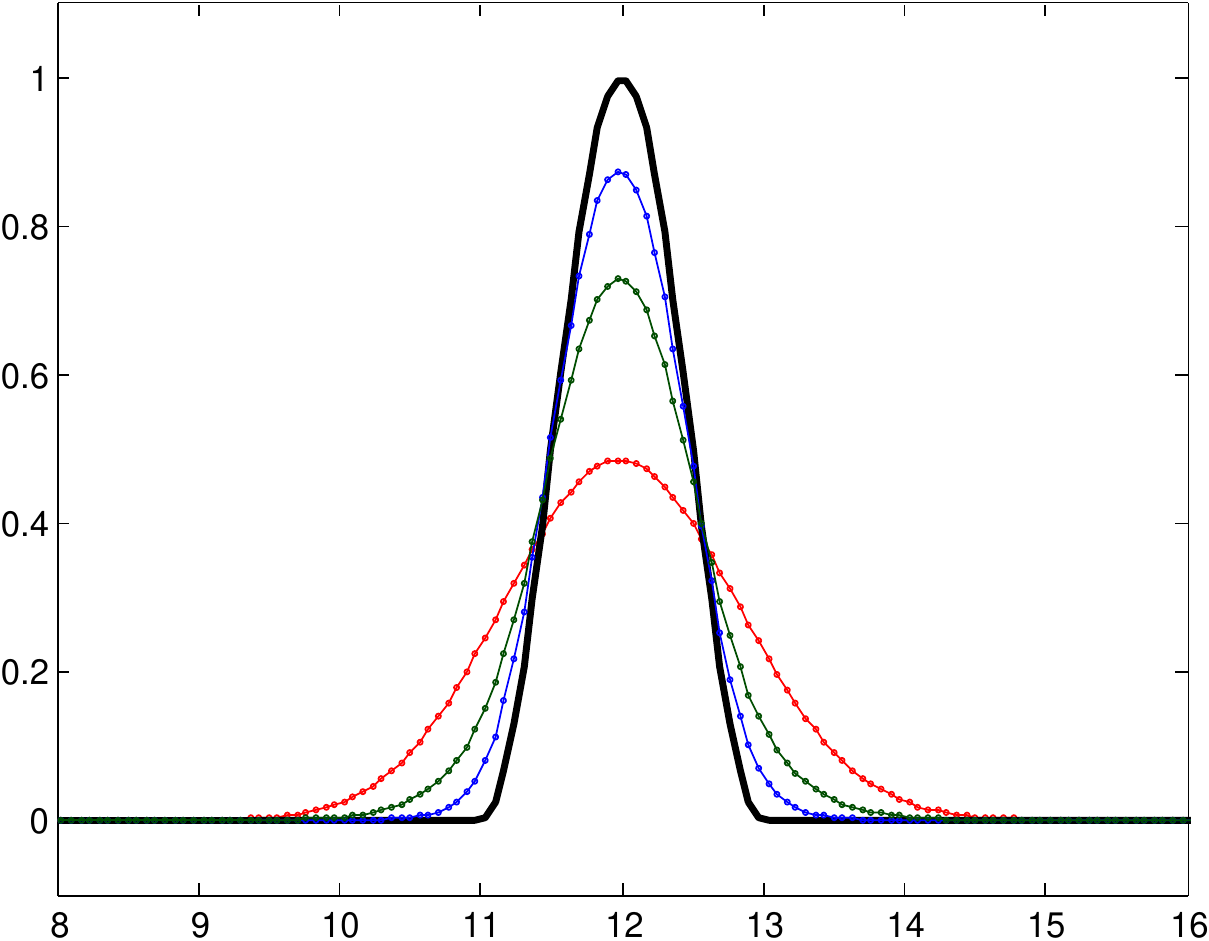}
\put(62,63){\footnotesize $\thetauno=0.99$}
\put(62,55){\footnotesize $\thetauno=0.997$}
\put(70,27){\footnotesize $\thetauno=1$}
\put(62,72){\footnotesize $u(\cdot,T)$}
\put(62,72){\line(-1,-1){8}}
\put(62,63){\line(-1,-1){7.5}}
\put(62,55){\line(-1,-1){7.2}}
\put(70,27){\line(-1,-2){4}}
\end{overpic} 
&
\begin{overpic}
[width=0.32\textwidth]{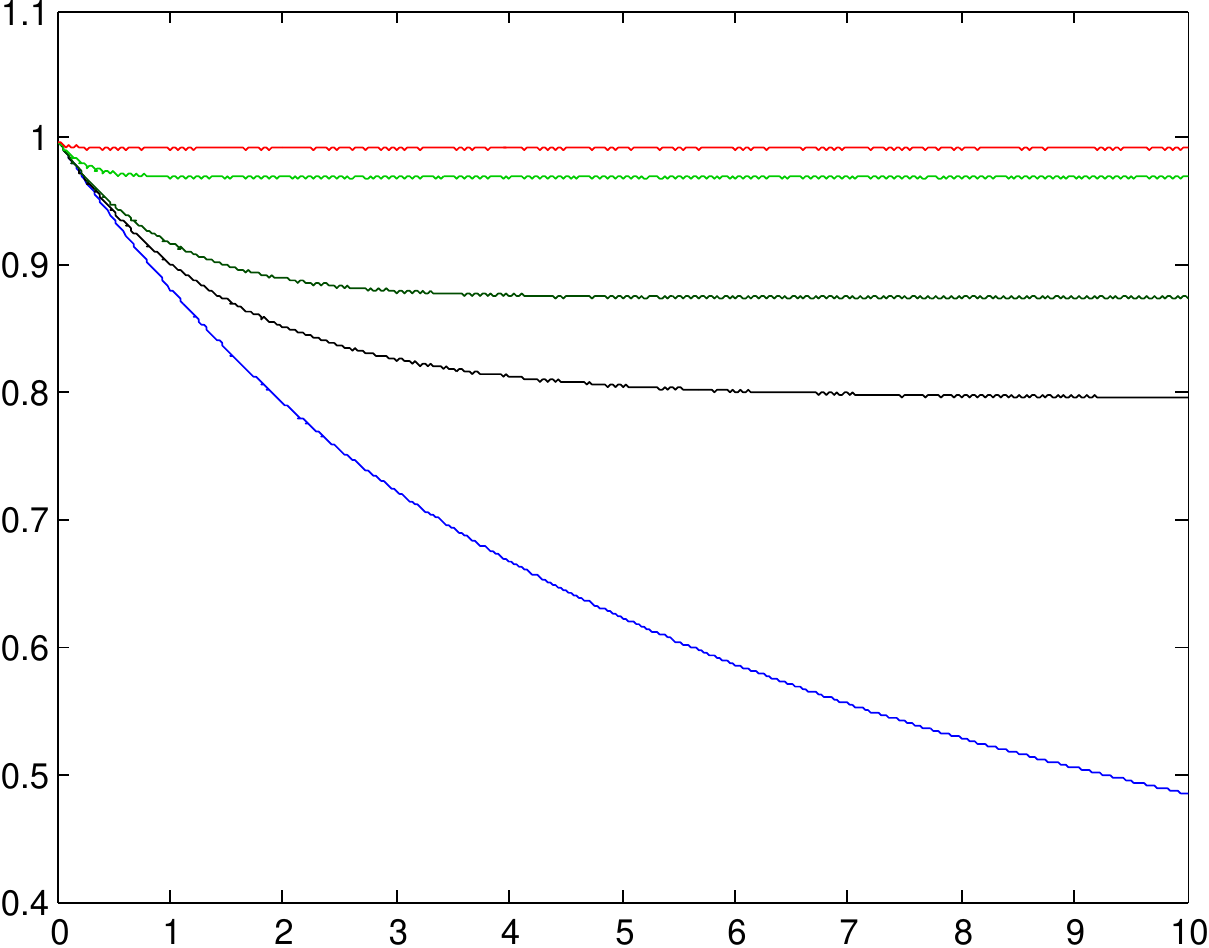}
\put(75,20){\footnotesize $\thetauno=1$}
\put(75,40){\footnotesize $\thetauno=0.995$}
\put(75,49){\footnotesize $\thetauno=0.99$}
\put(75,59){\footnotesize $\thetauno=0.95$}
\put(75,67){\footnotesize $\thetauno=0.8$}
\end{overpic} 
&
\begin{overpic}
[width=0.32\textwidth]{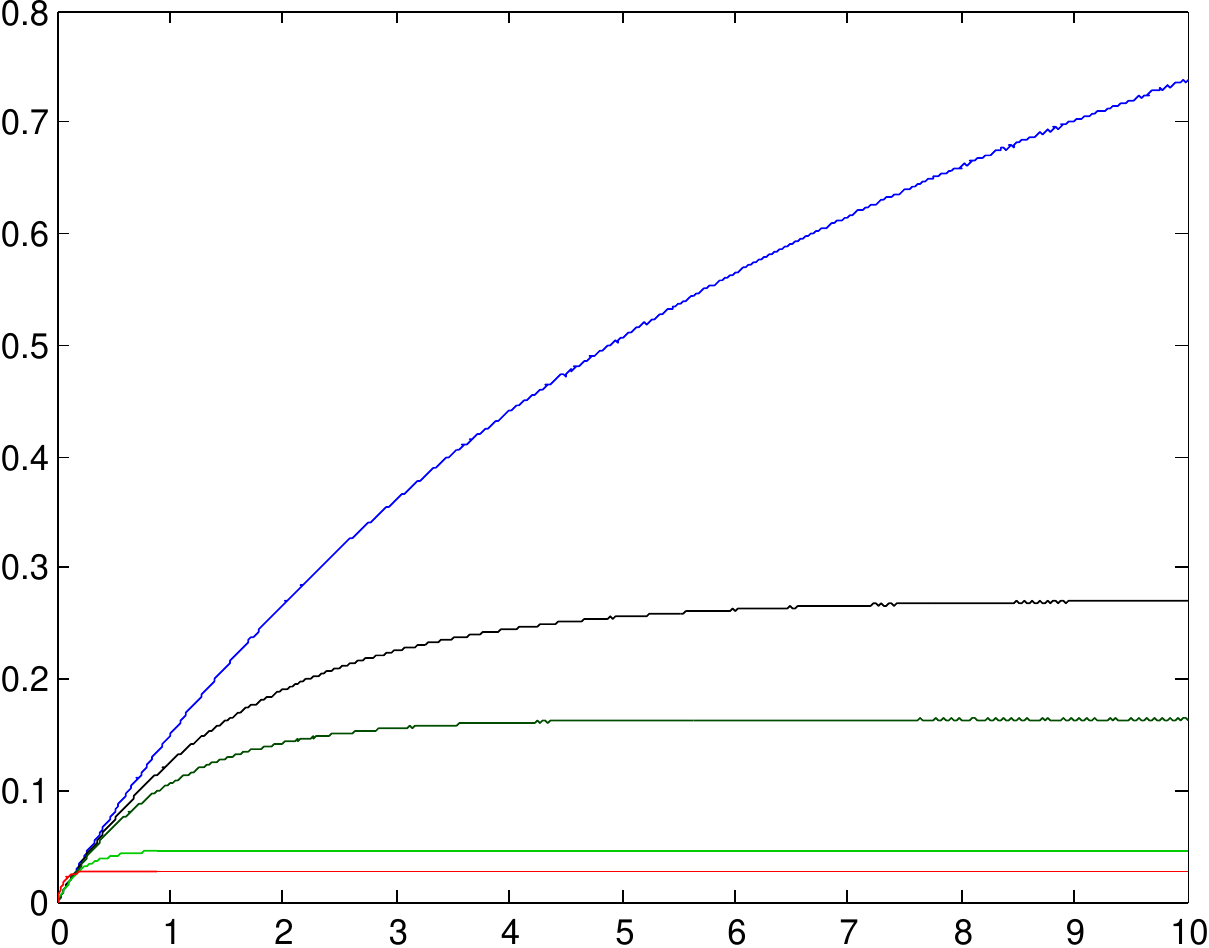}
\put(75,69){\footnotesize $\thetauno=1$}
\put(73,30){\footnotesize $\thetauno=0.995$}
\put(73,20){\footnotesize $\thetauno=0.99$}
\put(73,9){\footnotesize $\thetauno=0.95$}
\put(45,7){\footnotesize $\thetauno=0.8$}
\end{overpic}  
\\
(a) & (b) & (c)
\end{tabular}
\end{center}
\caption{Toy test II: (a) the exact solution $u$ and the approximate solution $W$ 
at final time, for $\thetadue=1$ and $\thetauno=0.997,0.99,1$. (b) max and (c) $L^1$ error of the approximate solution $W$ as a function of time, for $\thetadue=1$ and $\thetauno=0.8,0.95,0.99,0.995,1$.}
\label{fig:es1}
\end{figure}
Here the Courant number is equal to 0.1875. 
Fig.\ \ref{fig:es1}(a) shows the approximate solution $W$ at final time for different values of $\thetauno$. For $\thetauno=1$ (pure UPW) it is visible the well-known diffusive behaviour of the scheme. Switching on the coupling, namely assuming $\thetauno<1$, the behaviour changes completely. After a short transient during which the solution is smeared out, the solution is frozen and propagates at correct speed forever, with no further change in the shape. 
Fig.\ \ref{fig:es1}(b,c) shows the evolution of the maximum of the solution and the $L^1$ error for different values of $\thetauno$. For $\thetauno=1$, the maximum decreases, as a consequence of the numerical diffusion. For $\thetauno<1$, the maximum decreases for a short time, and then it comes to a regime state. Analogously, the $L^1$ error is bounded.

%%%%%%%%%%%%%%%%%%%%%%%%%%%%%%%%%%%%%%%
%%%%%%%%%%%%%%%%%%%%%%%%%%%%%%%%%%%%%%%
%%%%%%%%%%%%%%%%%%%%%%%%%%%%%%%%%%%%%%%
\section{Numerical tests}\label{sec:num_tests}
In this section we present some numerical tests, showing the behavior of the blended scheme \eqref{BlendedScheme} and confirming the theoretical results described in the previous sections. 
We always assume that the CFL condition applies, i.e.\ the time step is chosen in such a way that
\begin{equation}\label{CFL}
\Dt=\CFL \frac{\Dx}{\|A\|_\infty} \quad\text{(eq.\ \eqref{AE})}, \qquad\qquad
\Dt=\CFL \frac{\Dx}{\|f_u\|_\infty} \quad\text{(eq.\ \eqref{CL})},
\end{equation}
where $\CFL<1$ is chosen small enough to guarantee the stability of the scheme in use. 
Recalling definition \eqref{def:E1}, in the case of Eulerian-Eulerian coupling we define the \emph{reference error} $\Eref$ as
\begin{equation}\label{def:Eref}
\Eref:=\min\{E^1[W(1,1)],E^1[V(1,1)]\}\,,
\end{equation}
i.e.\ the minimum error one can achieve with the original uncoupled schemes. 
We also define the following subset of the parameter space
\begin{equation}\label{def:YES}
\begin{split}
\YES[W]:=\{(\thetauno,\thetadue)\in[0,1]^2~:~E^1[W(\thetauno,\thetadue)]<\Eref\}, \\
\YES[V]:=\{(\thetauno,\thetadue)\in[0,1]^2~:~E^1[V(\thetauno,\thetadue)]<\Eref\},
\end{split}
\end{equation}
i.e.\ the couples $(\thetauno,\thetadue)$ which improve the approximation with respect to both the uncoupled schemes.

In the case of multiscale coupling we always set $\mu=1$ (partial coupling) since it seems more convenient, both from the accuracy and CPU time point of view. 
Moreover, to be fair we define
\begin{equation}\label{def:YESmM}
\Eref:=E^1[W(1,1)] \qquad \text{ and } \qquad
\YES[W]:=\{\thetauno\in[0,1]~:~E^1[W(\thetauno,1)]<\Eref\},
\end{equation}
using only the Eulerian density $W(1,1)$ and not $V(1,1)$. This choice is motivated by the fact that the reconstruction of the Lagrangian density proposed in \eqref{S2micro}, despite it is a natural choice in the framework of the blended scheme, it is not in general the best approximation one can achieve (see, e.g., the SPH method). As a consequence, it is quite easy to get more accurate approximation of the Lagrangian density. 

In order to use Richardson extrapolation, we introduce a scale parameter $s\in(0,\frac12]$ to coarsen the grid in a proper way. Given parameters $N_C$, $N_T$ for a reference grid $\GG$, 
we define (denoting by $[\cdot]$ the upper integer part) the grid $\GG'$ such that $N_C'=[sN_C]$, $N_T'=[sN_T]$, 
and the grid $\GG''$ such that $N_C''=2[sN_C]=2N_C'$, $N_T''=2[sN_T]=2N_T'$. We also define $\Np'=[s\Np]$ and $\Np''=2[s\Np]=2\Np'$ to scale the number of particles, see section \ref{sec:theory}.

Hereafter we use the following acronyms:
UPW = Upwind (I order), RLW = Richtmyer two-step Lax-Wendroff (II order) \cite[Sect.\ 12.2]{levequebook}, WENO2 = Weighted Essentially Non-Oscillatory with linear interpolation and TVD Runge-Kutta 3 approximation in time (small constant $\varepsilon=10^{-9}$) (III order) \cite{liu1994JCP}, EE = Explicit Forward Euler (I order).

%%%%%%%%%%%%%%%%%%%%%%%%%%%%%%%%%%%%%%%%%%%%%%%%%%%
\subsection{Test \MMx}
%%%%%%%%%%%%%%%%%%%%%%%%%%%%%%%%%%%%%%%%%%%%%%%%%%%
\begin{center}
\begin{tabular}{|c|c|c|c|c|c|c|c|c|c|}
\hline
Test & $\S_1$ & $\S_2$ & $\bar u(x)$ & $\Omega$ & $T$ & $A(x)$ & $N_C$ & $N_T$ & $\CFL$
\\ \hline
\MMx & RLW & UPW & $\chi_{\left[\frac12,\frac32\right]}(x)$ & $[0,20]$ & 2.3 & $x$ & 1200 & 3000 & 0.92
\\ \hline 
\end{tabular}
\end{center}
\vskip0.4cm

\begin{figure}[h!]
\begin{center}
\begin{tabular}{cc}
\begin{overpic}
[width=0.4\textwidth]{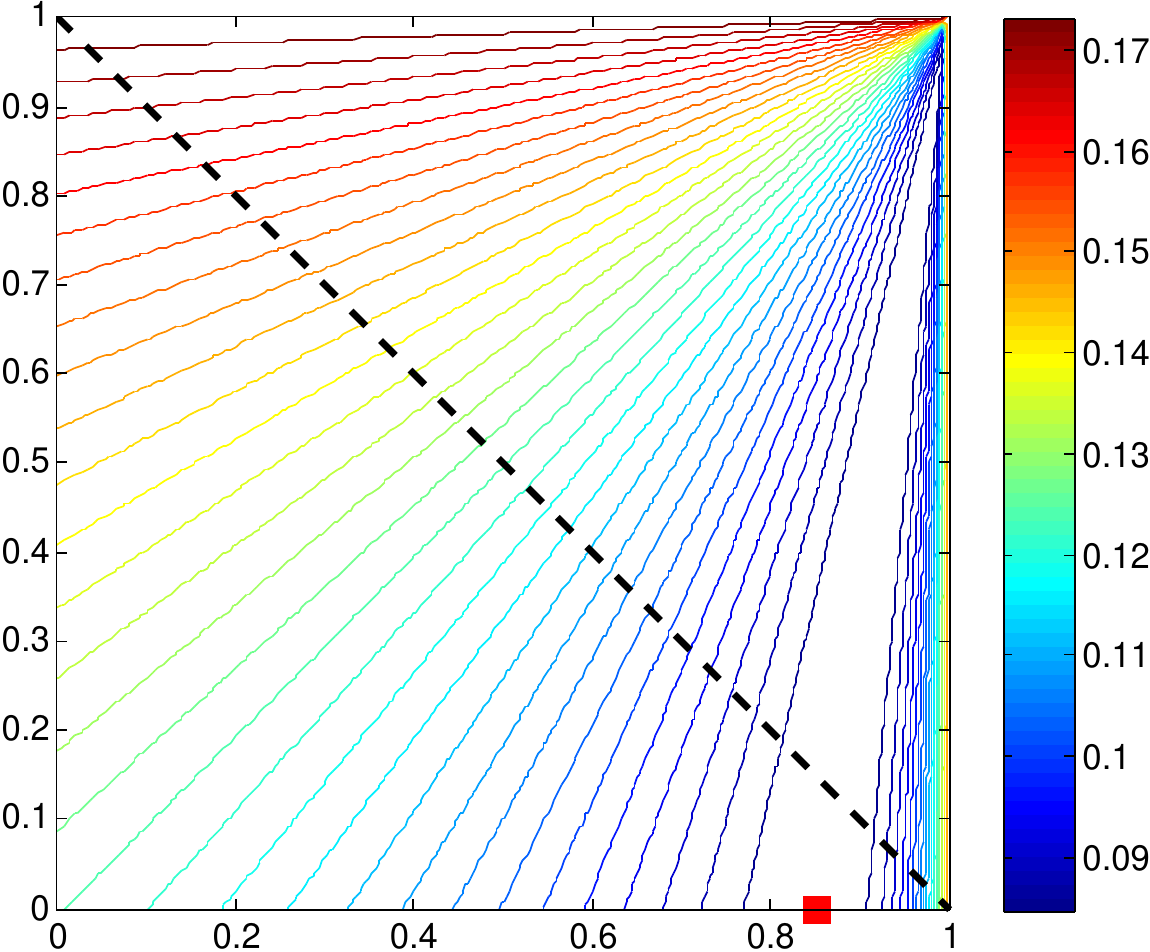}
\put(80,-4){$\thetauno$} \put(-4,79){$\thetadue$}
\end{overpic} 
\qquad\qquad &
\begin{overpic}
[width=0.331\textwidth]{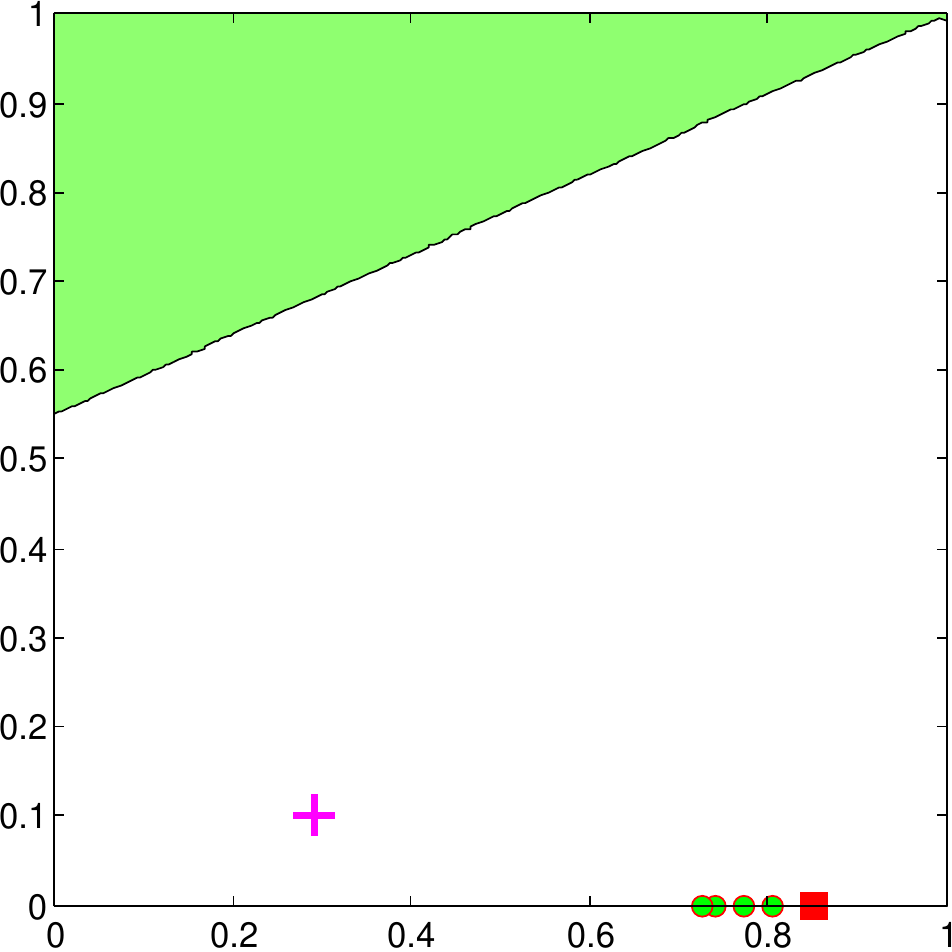}
\put(95,-5){$\thetauno$} \put(-5,95){$\thetadue$} \put(50,50){$\YES[W]$}
\end{overpic}  
\\
(a) & (b)
\end{tabular}
\end{center}
\caption{Test \MMx: (a) level sets of the function $(\thetauno,\thetadue)\to E^1[W(\thetauno,\thetadue)]$. (b) Region $\YES[W]$ (white). The red square indicates the global minimum point $(\thetauno^*,\thetadue^*)$, the green circles indicates the global minimum points $(\thetauno^*,\thetadue^*)$ on coarser grids ($s=\frac12,\frac13,\frac14,\frac15$), and the magenta cross indicates the Richardson minimum point $(\thetauno^R,\thetadue^R)$.}
\label{fig:MMx:batch-LS-YES}
\end{figure}

In this test we blend RLW and UPW in case of the positive space-dependent velocity $A(x)=x$. The exact solution is given by $u(x,t)=\bar u(xe^{-t})e^{-t}$. 
Fig.\ \ref{fig:MMx:batch-LS-YES} shows (a) the level sets of the function 
$(\thetauno,\thetadue)\to E^1[W(\thetauno,\thetadue)]$ and (b) the region $\YES[W]$. 
In this case one can note that the level sets of the error are straight lines and one can obtain the same solution varying the coupling parameters while keeping their ratio fixed. 
This implies that the line $\lambda=1-\mu$ nearly cuts all the levels sets, and, therefore, restriction to such a particular case (much simpler to study analytically since the system is reduced to a single equation) does not lead to a loss of generality. 

The reference error $\Eref$ in this case is equal to 0.1463, while the minimum error that can be achieved is 0.0816 (-44.22\%), corresponding to the couple $(\thetauno^*,\thetadue^*)=(0.8533,0)$ computed by an exhaustive search in the parameter space. Richardson extrapolation method (with $s=\frac18$) suggests instead the couple $(\thetauno^R,\thetadue^R)=(0.29,0.1)$, which is fully inside $\YES[W]$ and leads to the error 0.117 (-20.03\%). Fig.\ \ref{fig:MMx:batch-LS-YES}(b) also shows the best coupling parameters $(\thetauno^*,\thetadue^*)$ computed coarsening the grid with $s=\frac12,\frac13,\frac14,\frac15$. It can be seen that the best coupling parameters do not depend very much on the grid and move monotonically with respect to the scaling from right to left.

\begin{figure}[h!]
\begin{center}
\begin{tabular}{ccc}
\includegraphics[width=0.32\textwidth]{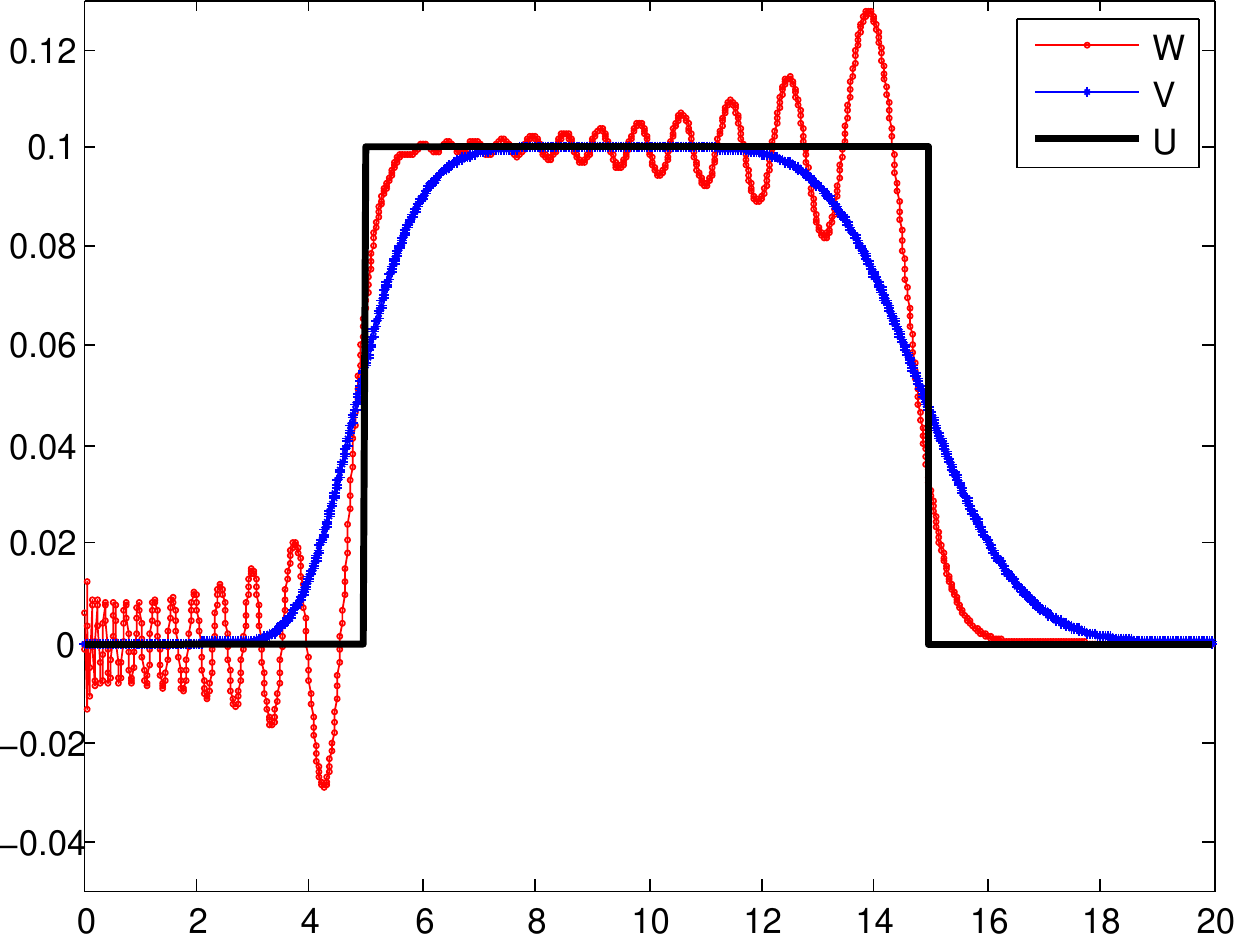} &
\includegraphics[width=0.32\textwidth]{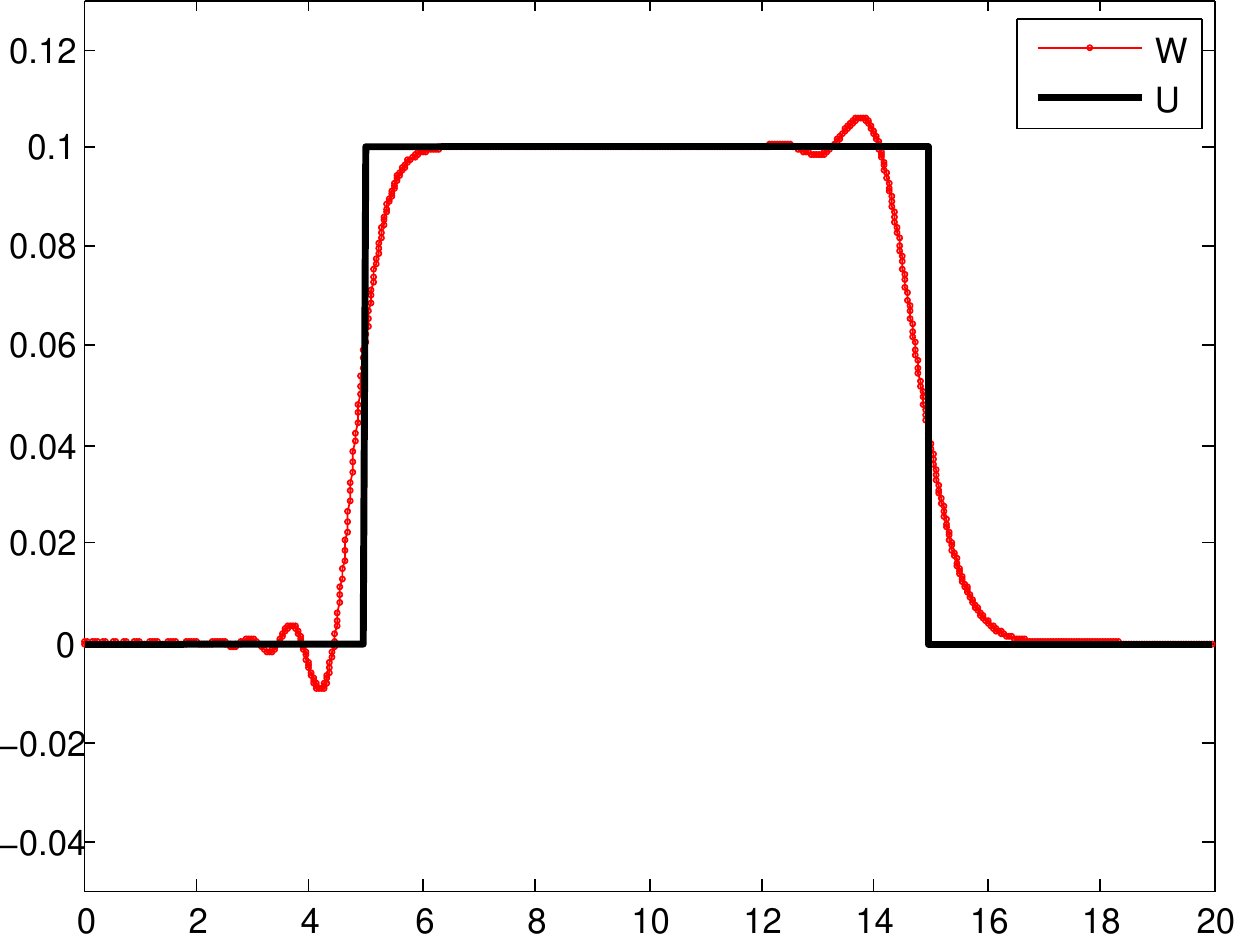} &
\includegraphics[width=0.32\textwidth]{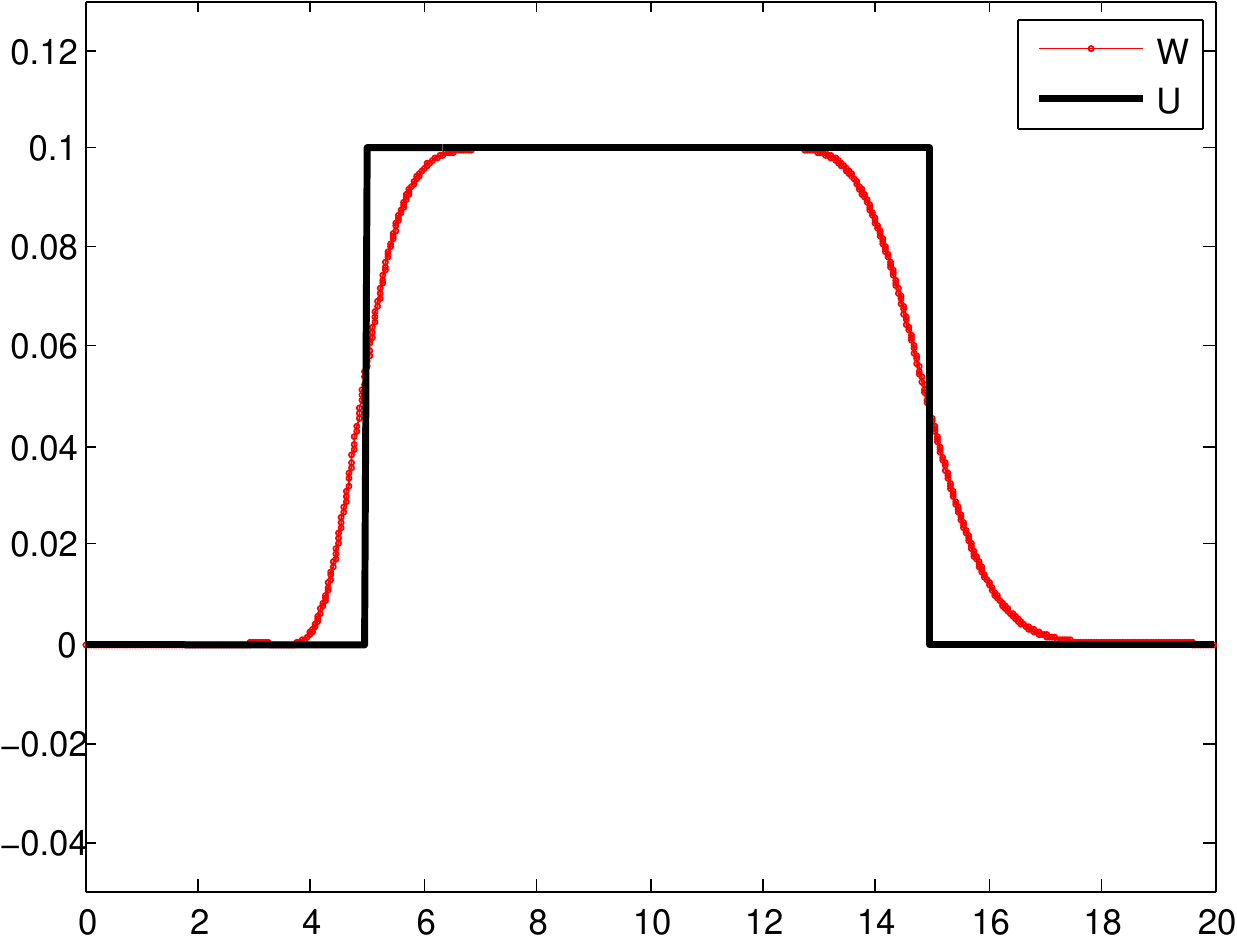} \\
(a) & (b) & (c)
\end{tabular}
\end{center}
\caption{Test \MMx: comparison with the exact solution. 
(a) $W^{N_T}(1,1)$ and $V^{N_T}(1,1)$. 
(b) $W^{N_T}(\thetauno^*,\thetadue^*)$. 
(c) $W^{N_T}(\thetauno^R,\thetadue^R)$.}
\label{fig:MMx:pures_and_bests}
\end{figure}

Fig.\ \ref{fig:MMx:pures_and_bests} shows the exact solution $U^{N_T}$ together with (a) the two uncoupled solutions $W^{N_T}(1,1)$ and $V^{N_T}(1,1)$, (b) the best coupling $W^{N_T}(\thetauno^*,\thetadue^*)$, and (c) the best Richardson coupling $W^{N_T}(\thetauno^R,\thetadue^R)$ at final time. 
Note that, although $W^{N_T}(\thetauno^R,\thetadue^R)$ is worst than $W^{N_T}(\thetauno^*,\thetadue^*)$ in terms of $L^1$ error, it does not show spurious oscillations. This is due to the fact that oscillations are sensitive to the grid and then Richardson method easily discards such a solutions. 

To better quantify the improvement of the blended scheme with respect to the uncoupled schemes, we found the grid needed by the $L^1$-best uncoupled scheme (RLW) to achieve the same accuracy of $W^{N_T}(\thetauno^R,\thetadue^R)$ and $W^{N_T}(\thetauno^*,\thetadue^*)$. In the first case we need $N_C=1704$, $N_T=4260$, while in the second case $N_C=3240$, $N_T=8100$, corresponding respectively to a refinement of the grid of a factor $1.42$ and $2.71$. 

%%%%%%%%%%%%%%%%%%%%%%%%%%%%%%%%%%%%%%%%%%%%%%%%%%%
\subsection{Test \mMx}
%%%%%%%%%%%%%%%%%%%%%%%%%%%%%%%%%%%%%%%%%%%%%%%%%%%
\begin{center}
\begin{tabular}{|c|c|c|c|c|c|c|c|c|c|c|}
\hline
Test & $\S_1$ & $\S_2$ & $\bar u(x)$ & $\Omega$ & $T$ & $A(x)$ & $N_C$ & $N_T$ & $\Np$ & $\CFL$
\\ \hline
\mMx & UPW & EE & $\chi_{\left[\frac12,\frac32\right]}(x)$ & $[0,20]$ & 2.3 & $x$ & 1200 & 3000 & $5\times N_C$ & 0.92
\\ \hline 
\end{tabular}
\end{center}
\vskip0.4cm

\begin{figure}[h!]
\begin{center}
\begin{tabular}{cc}
\begin{overpic}
[width=0.45\textwidth]{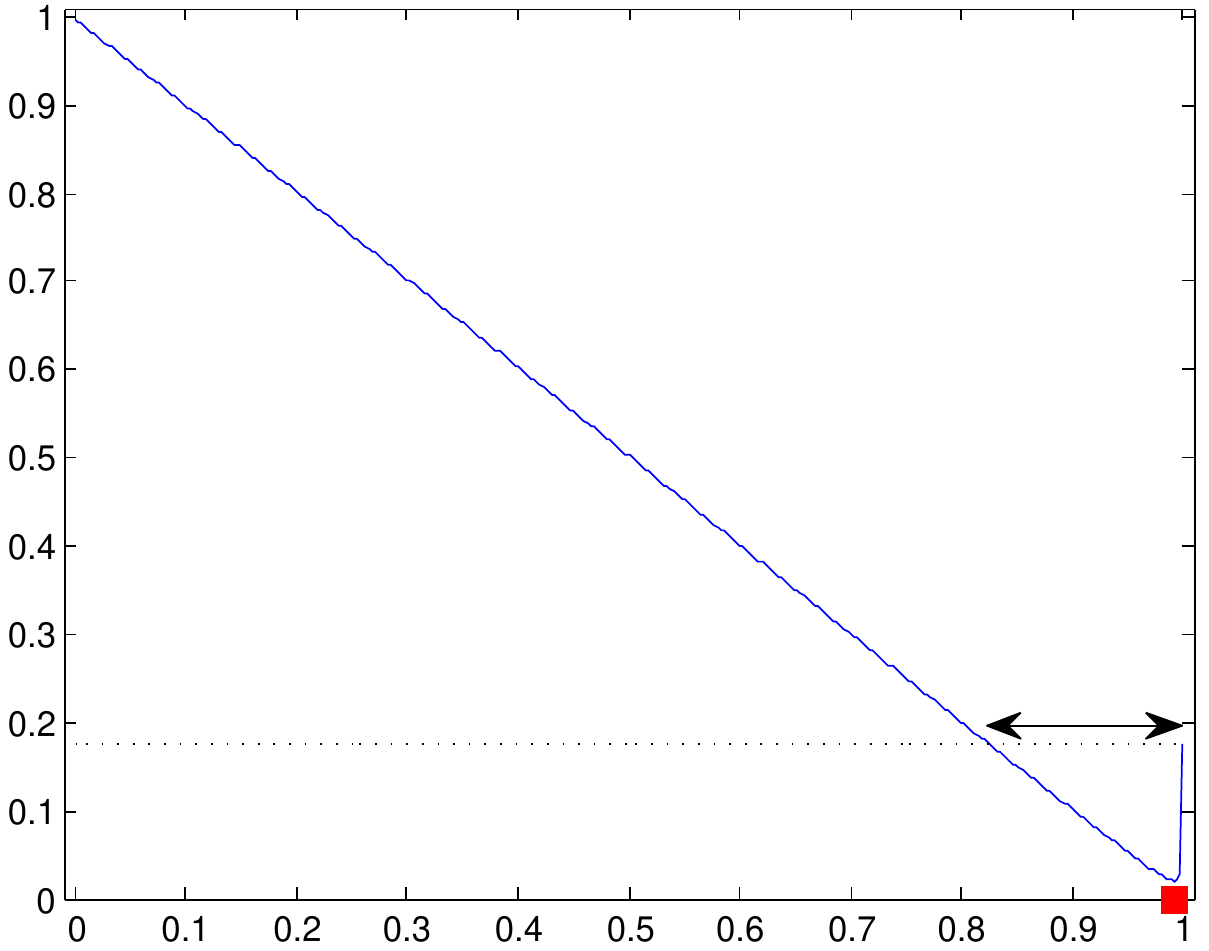}
\put(85,22){$\YES[W]$} 
\end{overpic}  
&
\includegraphics[width=0.45\textwidth]{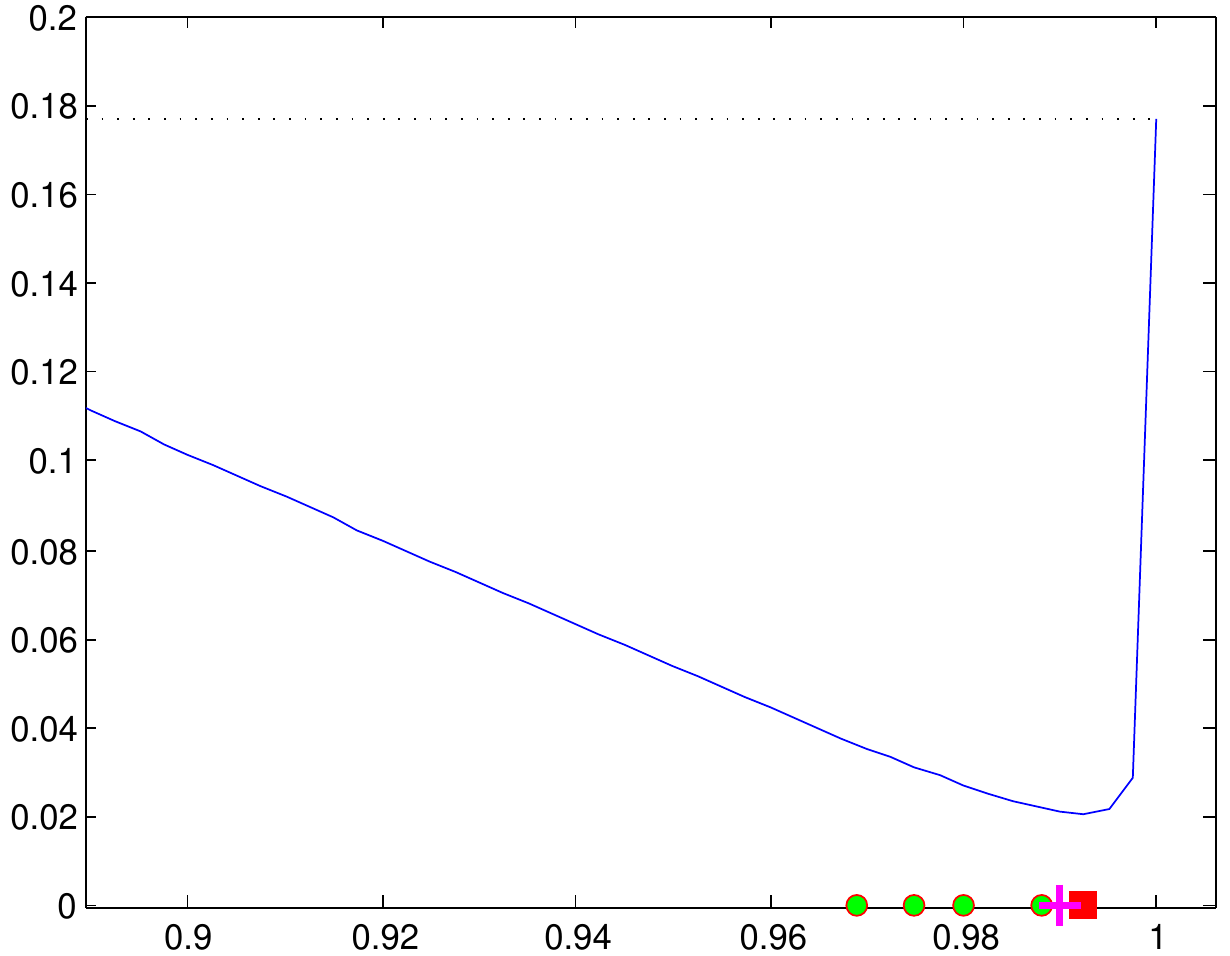} \\
(a) & (b)
\end{tabular}
\end{center}
\caption{Test \mMx: function $\thetauno \to E^1[W(\thetauno,1)]$ cut by the reference error. (a) $\thetauno\in[0.6,1]$. (b) $\thetauno\in[0.89,1]$. The red square indicates the global minimum point $\thetauno^*$, the green circles indicates the global minimum points $\thetauno^*$ on coarser grids ($s=\frac12,\frac14,\frac16,\frac18$), and the magenta cross indicates the Richardson minimum point $\thetauno^R$.}
\label{fig:mMx:batch}
\end{figure}

In this test we run the multiscale version of the blended scheme coupling UPW and EE in the same setting of Test \MMx. 
The UPW error, which will be taken as the reference one, is equal to 0.1771.   
The minimum error that can be achieved is 0.0204 (-88.48\%), corresponding to $\thetauno^*=0.992$ computed by an exhaustive search in the parameter space $[0,1]\times\{1\}$. 
Fig.\ \ref{fig:mMx:batch} shows the function $\thetauno \to E^1[W(\thetauno,1)]$ cut by the reference error.
Richardson extrapolation method (with $s=\frac13$) suggests $\thetauno^R=0.99$ which leads to the error 0.0208 (-88.26\%).

\begin{figure}[h!]
\begin{center}
\begin{tabular}{cc}
\includegraphics[width=0.45\textwidth]{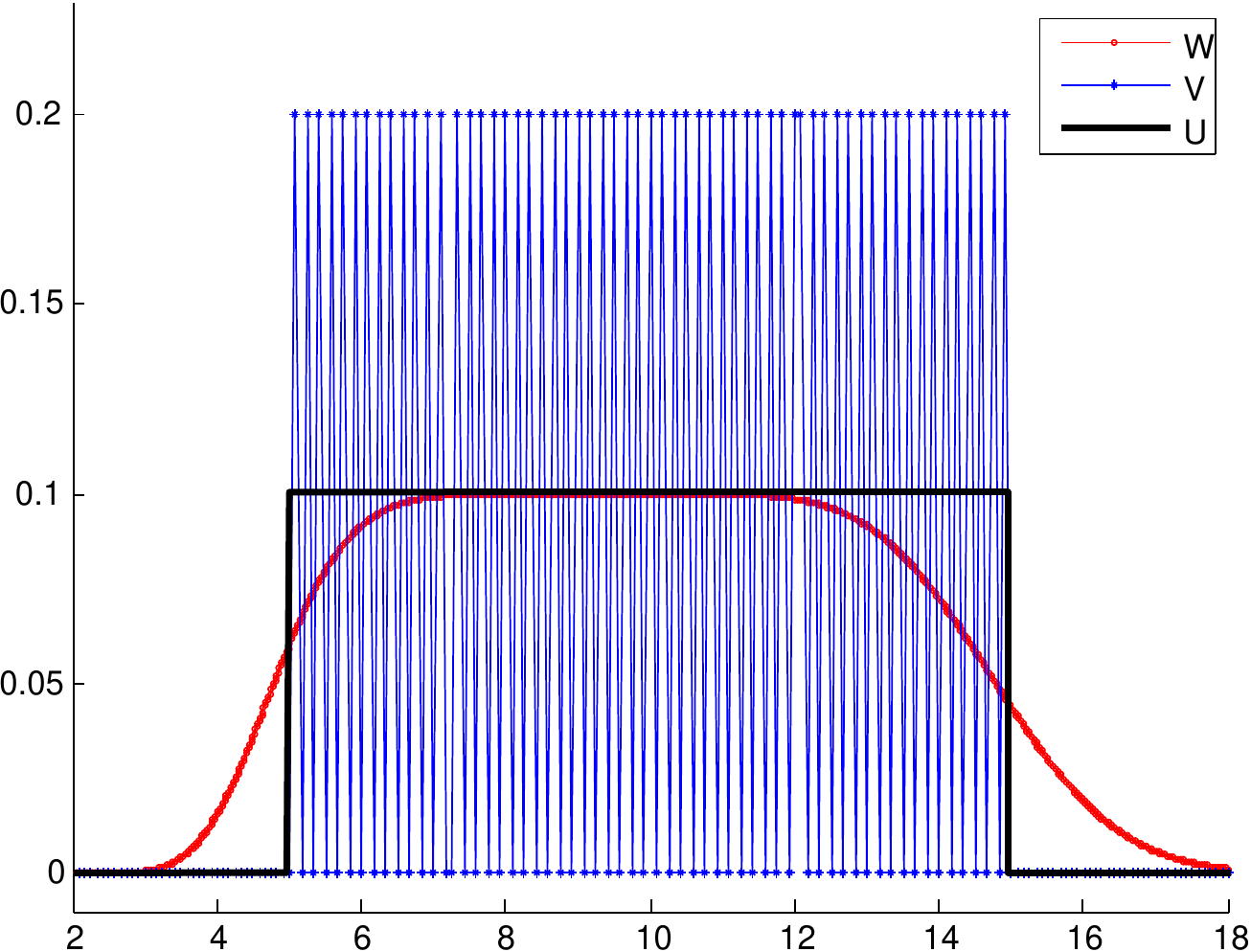} &
\includegraphics[width=0.45\textwidth]{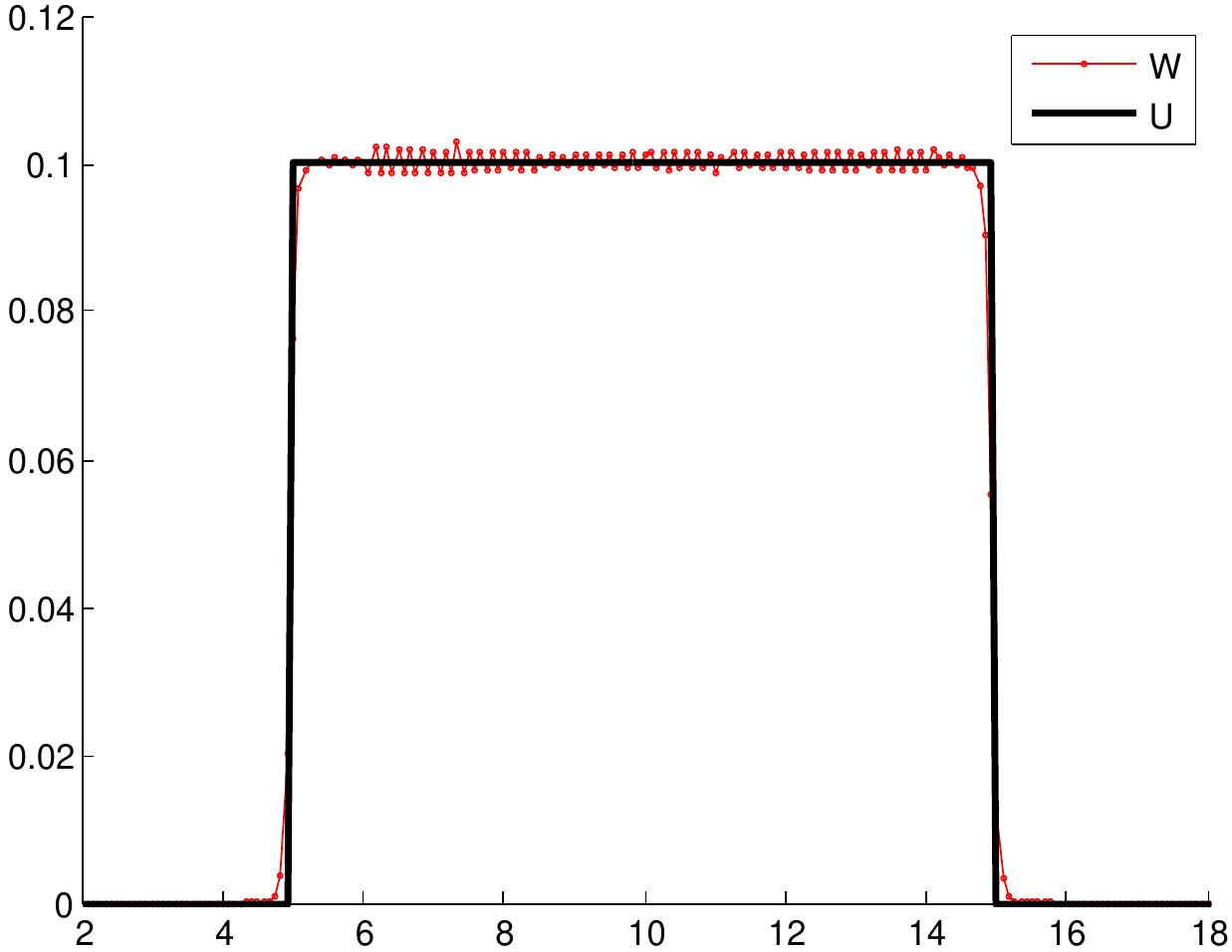} \\
(a) & (b)
\end{tabular}
\end{center}
\caption{Test \mMx: comparison with the exact solution. 
(a) $W^{N_T}(1,1)$ and $V^{N_T}(1,1)$. 
(b) $W^{N_T}(\thetauno^R,1)$. 
Solutions are downsampled for a better viewing.}
\label{fig:mMx:pures_and_bests}
\end{figure}

Fig.\ \ref{fig:mMx:pures_and_bests} shows the exact solution $U^{N_T}$ together with (a) the two uncoupled solutions $W^{N_T}(1,1)$ and $V^{N_T}(1,1)$, and (b) the best Richardson coupling $W^{N_T}(\thetauno^R,1)$ at final time. 
It can be seen that EE, as opposed to UPW, is rather good in terms of diffusion and in capturing discontinuities, but it has several oscillations. The size of the oscillations depends on $\Np$, while their number depends on $N_C$. The vector field $A(x)=x$ tends to separate the particles (but in this case they remains equispaced), therefore the oscillations become larger as the time goes by. Note the the Lagrangian density $V(1,1)$ is 0 in those cells not covered by any particle. The blended solution $W^{N_T}(\thetauno^R,1)$ combines nicely the advantages of the two schemes, showing little diffusion and tiny oscillations (which have a different nature with respect to those of high-order Eulerian schemes).
For the sake of comparison, we report that WENO2 scheme reaches the same accuracy of UPW+EE on a grid refined by a factor 4.2.

\begin{figure}[h!]
\begin{center}
\begin{tabular}{cc}
\includegraphics[width=0.45\textwidth]{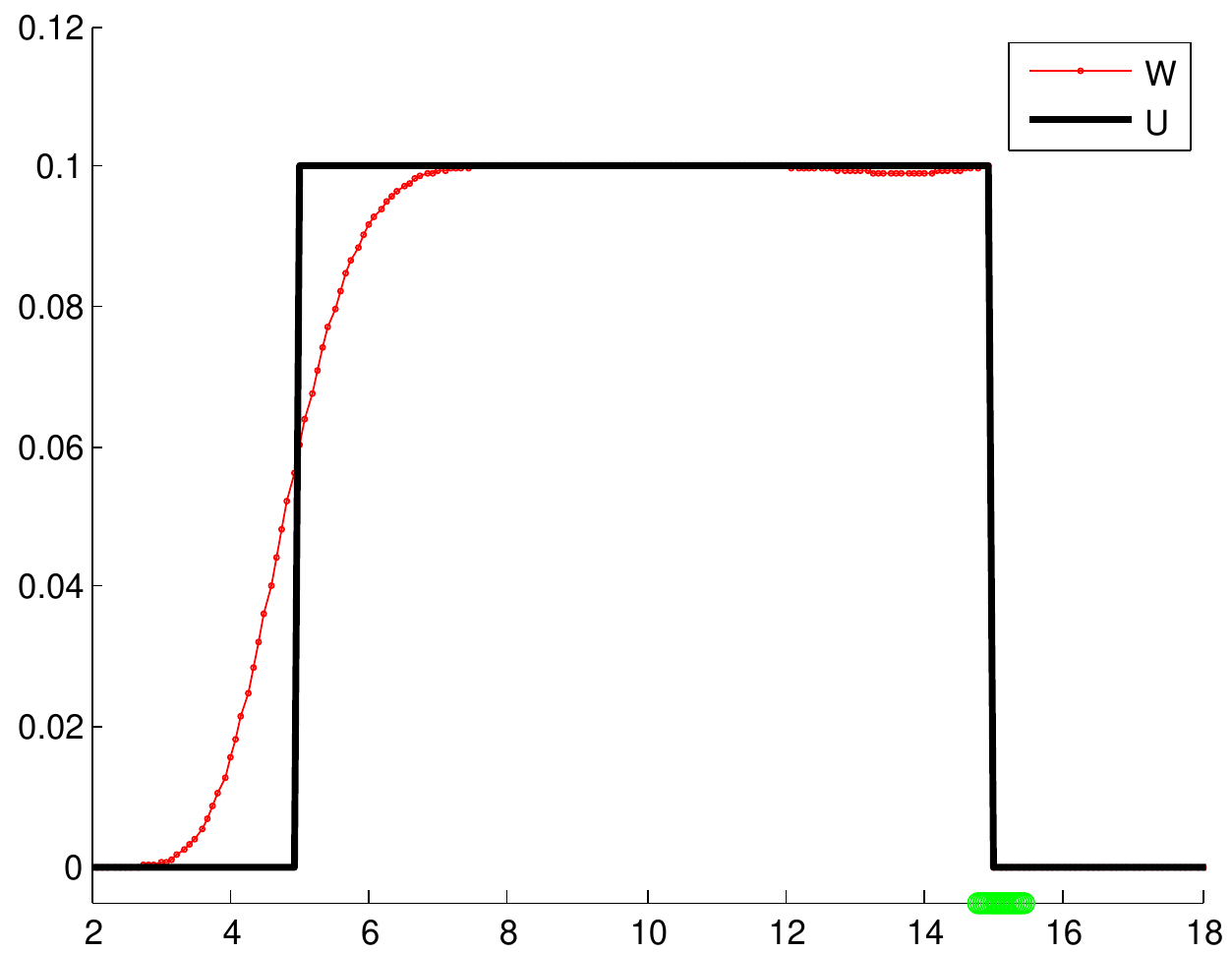} &
\includegraphics[width=0.45\textwidth]{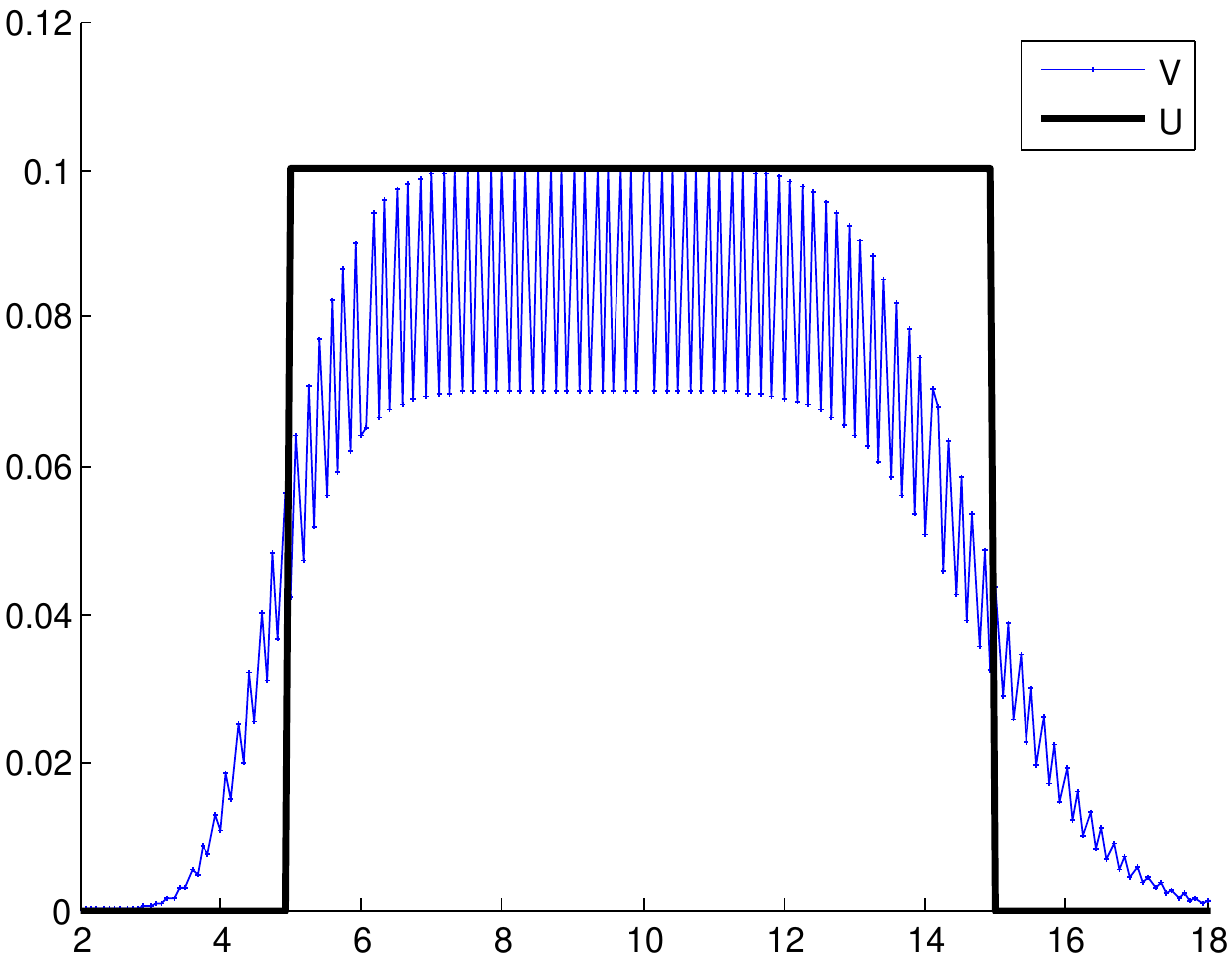} \\
(a) & (b)
\end{tabular}
\end{center}
\caption{Test \mMx: (a) $W^{N_T}(0,1)$ with localized particles. Green circles are the particles. (b) $V^{N_T}(1,0.3)$.}
\label{fig:mMx:localized}
\end{figure}

We also tested the possibility of localizing the particles around critical zones, namely the discontinuities. At initial time, we locate a few particles across the right discontinuity, and then we activate the blend (with $\thetauno=0$) only in the cells which are between the first and the last particle (more in general, inside the support of the Lagrangian density). This is needed in order to distinguish the two kind of no-particle zones, the one due to the particle rarefaction and the one uncovered by user's choice. 
In Fig.\ \ref{fig:mMx:localized}(a) we show the result of the blended scheme in this case and the region actually covered by particles. The right discontinuity is well caught despite the small number of particles.

For the sake of completeness, we also report the outcome of the blended scheme whenever the Lagrangian solution is regularized by the Eulerian one. In Fig.\ \ref{fig:mMx:localized}(b) we show the solution $V$ obtained with $(\thetauno,\thetadue)=(1,0.3)$. In this case the blend puts together, rather than cancels, the drawbacks of the two uncoupled schemes, leading to a solution which is both diffusive and oscillating. 

\newpage
%%%%%%%%%%%%%%%%%%%%%%%%%%%%%%%%%%%%%%%%%%%%%%%%%%%
\subsection{Test \mMsinx}
%%%%%%%%%%%%%%%%%%%%%%%%%%%%%%%%%%%%%%%%%%%%%%%%%%%
\begin{center}
\begin{tabular}{|c|c|c|c|c|c|c|c|c|c|c|}
\hline
Test & $\S_1$ & $\S_2$ & $\bar u(x)$ & $\Omega$ & $T$ & $A(x)$ & $N_C$ & $N_T$ & $\Np$ & $\CFL$
\\ \hline
\mMsinx & UPW & EE & $\sin(e^{2x}/20)$ & $[0,\pi]$ & $1$ & $\sin(x)$ & 600 & 200 & $1\times N_C$ & 0.96
\\ \hline 
\end{tabular}
\end{center}
\vskip0.4cm

\begin{figure}[h!]
\begin{center}
\begin{tabular}{cc}
\includegraphics[width=0.4\textwidth]{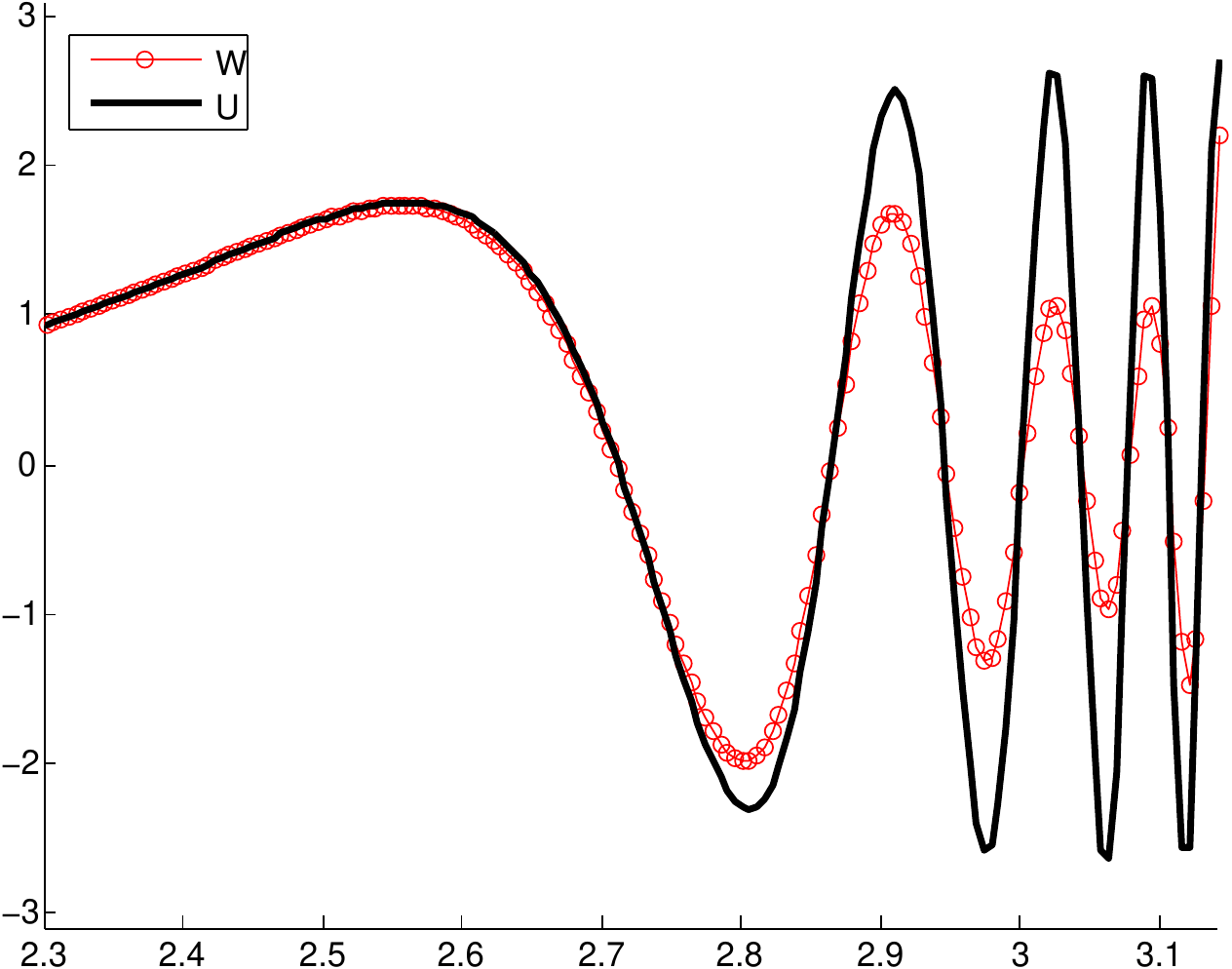} &
\includegraphics[width=0.4\textwidth]{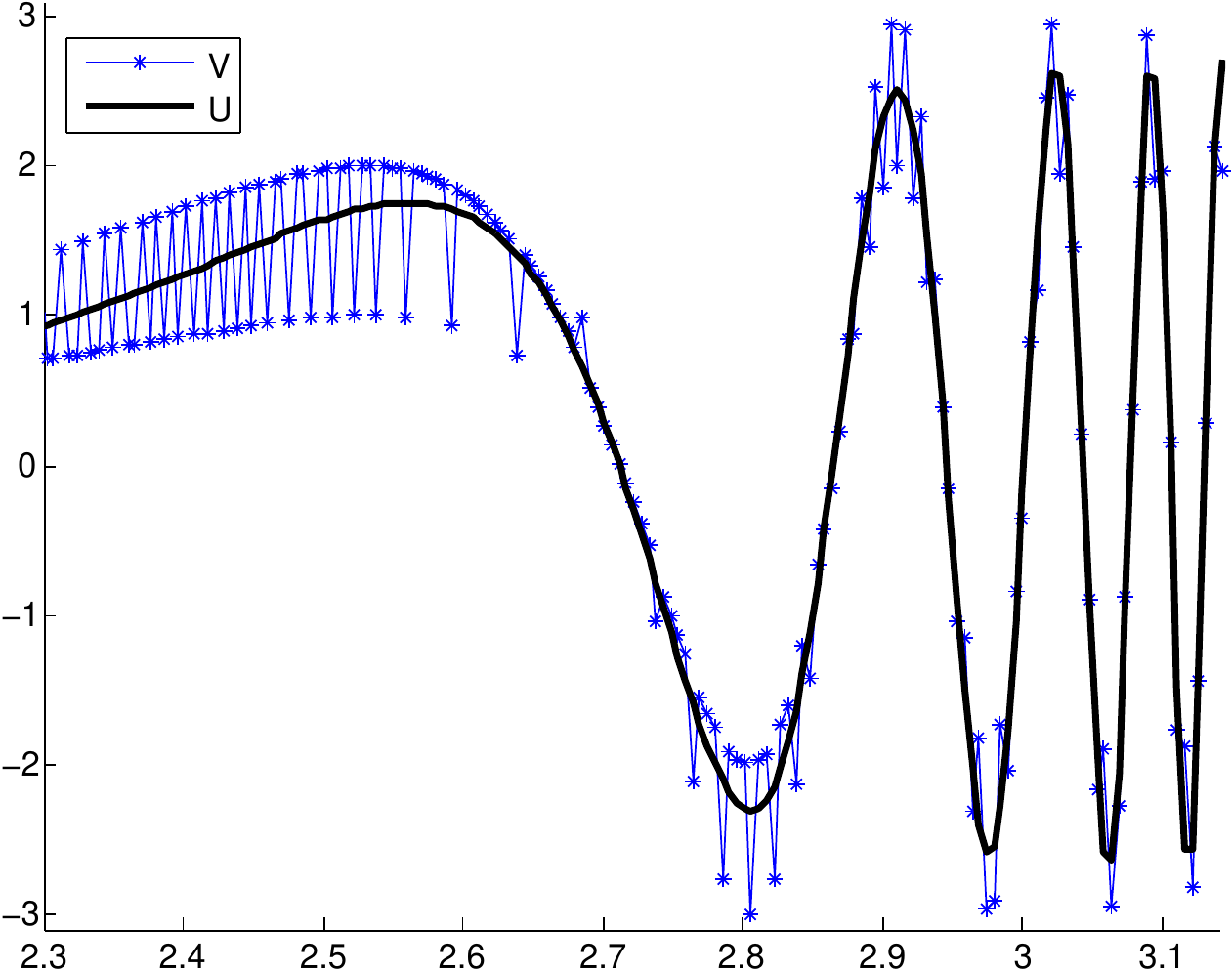} \\
(a) & (b) 
\end{tabular}
\end{center}
\caption{Test \mMsinx: uncoupled schemes in $[2.3,\pi]$. 
(a) $W^{N_T}(1,1)$. 
(b) $V^{N_T}(1,1)$.}
\label{fig:mMsinx:pures}
\end{figure}

In this test we run the multiscale version of the blended scheme coupling UPW and EE in case of a positive space-dependent velocity $A(x)=\sin(x)$ (for $x\in\Omega=[0,\pi]$). The exact solution is 
$$
u(x,t)=\bar u(\gamma(x,t))\left(\frac12\tan\left(\frac{x}{2}\right)\left(1+e^{-2t}\right)\sin(\gamma(x,t))+e^{-t}\cos(\gamma(x,t))\right),\quad
\gamma(x,t)=2\arctan\left(e^{-t}\tan\left(\frac{x}{2}\right)\right).
$$ 
The reference error (UPW) in this case is equal to 0.2591. 
Fig.\ \ref{fig:mMsinx:pures} shows the exact solution $U^{N_T}$ together with the two uncoupled solutions $W^{N_T}(1,1)$ and $V^{N_T}(1,1)$ at final time. In this case the particles move rightward accumulating on the right side of the domain. 
The minimum error that can be achieved is 0.0731 (-71.79\%), corresponding to $\thetauno^*=0.93$. Richardson extrapolation method (with $s=\frac12$) suggests instead $\thetauno^R=0.916$, which leads to the error 0.0742 (-71.36\%). 

\begin{figure}[b!]
\begin{center}
\includegraphics[width=0.98\textwidth]{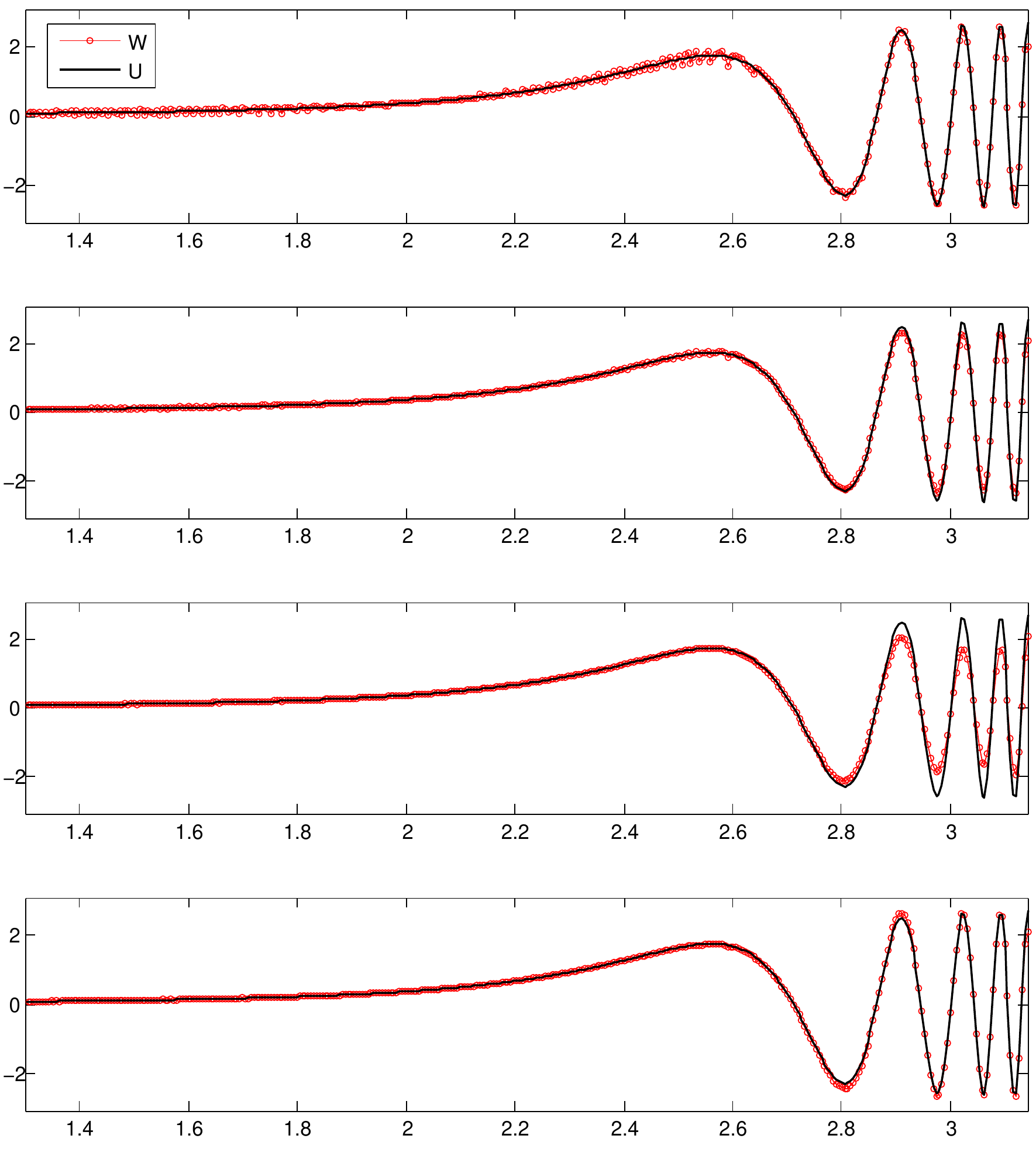}
\end{center}
\caption{Test \mMsinx: $W^{N_T}(\thetauno,1)$ in $[1.3,\pi]$ for different choices of $\thetauno$. from top to bottom: $\thetauno=0.8$ ($E^1=0.1085$), $\thetauno=\thetauno^R=0.952$ ($E^1=0.0764$), $\thetauno=0.99$ ($E^1=0.1517$), $\thetauno=\thetauno(\Lambda)$ ($E^1=0.0396$).}
\label{fig:mMsinx:4plots}
\end{figure}

Fig.\ \ref{fig:mMsinx:4plots} shows four blended solutions obtained with different choices of $\thetauno$. For low values of the coupling parameter ($\thetauno=0.8$) the solution is accurate where the number of particles is large (right side) but not where there are no or a few particles (left side). The situation is reversed for large values of the coupling parameter ($\thetauno=0.99$). The Richardson solution $\thetauno=\thetauno^R=0.916$ reaches a good compromise between the two cases. On the other hand, an even better result can be obtained by means of a variable-in-space $\thetauno=\thetauno(\Lambda)$ (see Section \ref{sec:multiscalecoupling} for the definition of $\Lambda$). Last plot in Fig.\ \ref{fig:mMsinx:4plots} shows the result of an hand-tuned decreasing function $\thetauno(\Lambda)$ which leads to a good accuracy in both sides of the domain.

\newpage
%%%%%%%%%%%%%%%%%%%%%%%%%%%%%%%%%%%%%%%%%%%%%%%%%%%
\subsection{Test \mMlwr}
%%%%%%%%%%%%%%%%%%%%%%%%%%%%%%%%%%%%%%%%%%%%%%%%%%%
\begin{center}
\begin{tabular}{|c|c|c|c|c|c|c|c|c|c|c|}
\hline
Test & $\S_1$ & $\S_2$ & $\bar u(x)$ & $\Omega$ & $T$ & $f(u)$ & $N_C$ & $N_T$ & $\Np$ & $\CFL$
\\ \hline
\mMlwr & GODUNOV & EE & $\frac12\chi_{\left[0,2\right]}(x)$ & $[-0.2,7]$ & $4$ & $u(1-u)$ & 100 & 200 & $5\times N_C$ & 0.28
\\ \hline 
\end{tabular}
\end{center}
\vskip0.4cm

\begin{figure}[h!]
\begin{center}
\begin{tabular}{cc}
\includegraphics[width=0.45\textwidth]{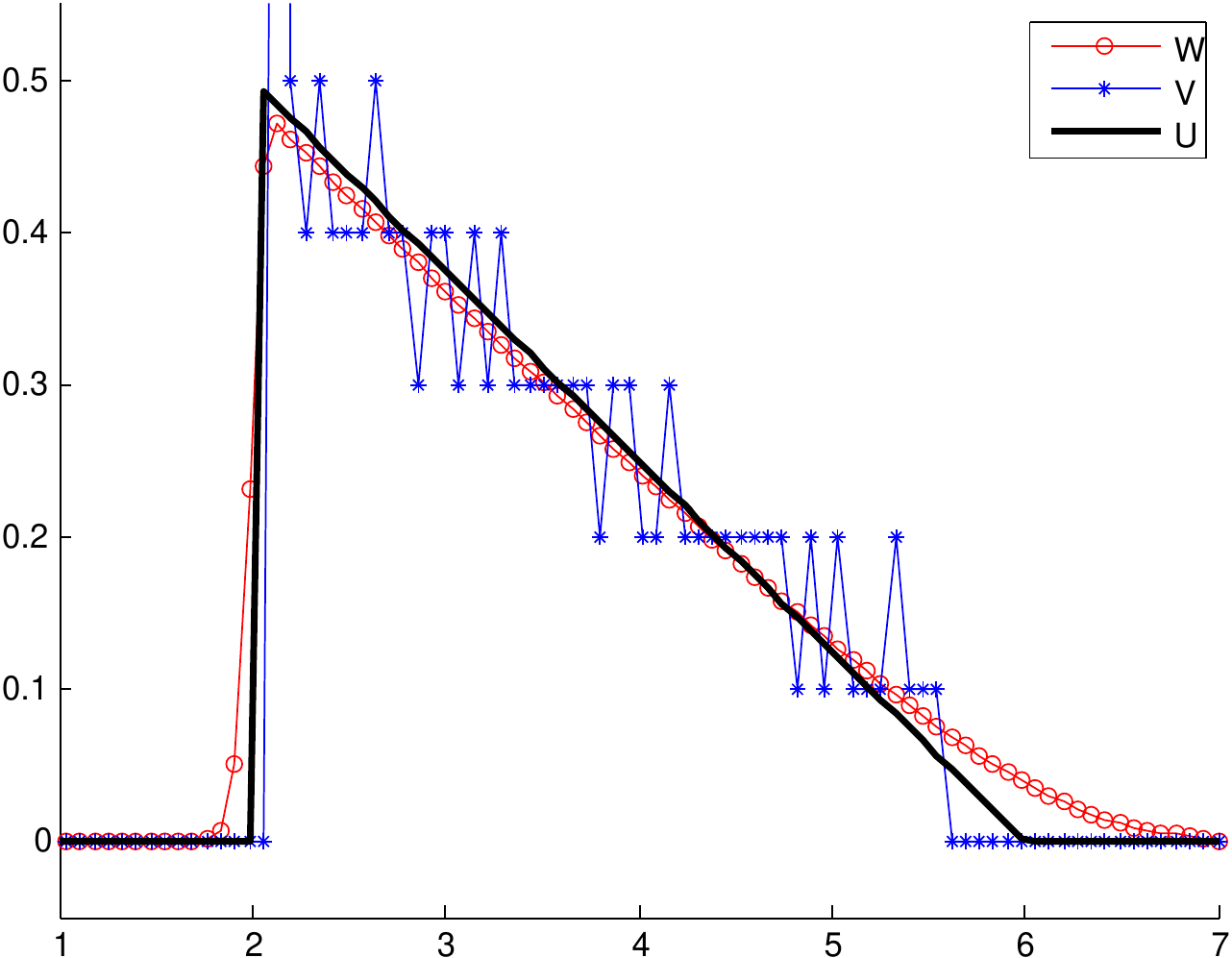} &
\includegraphics[width=0.45\textwidth]{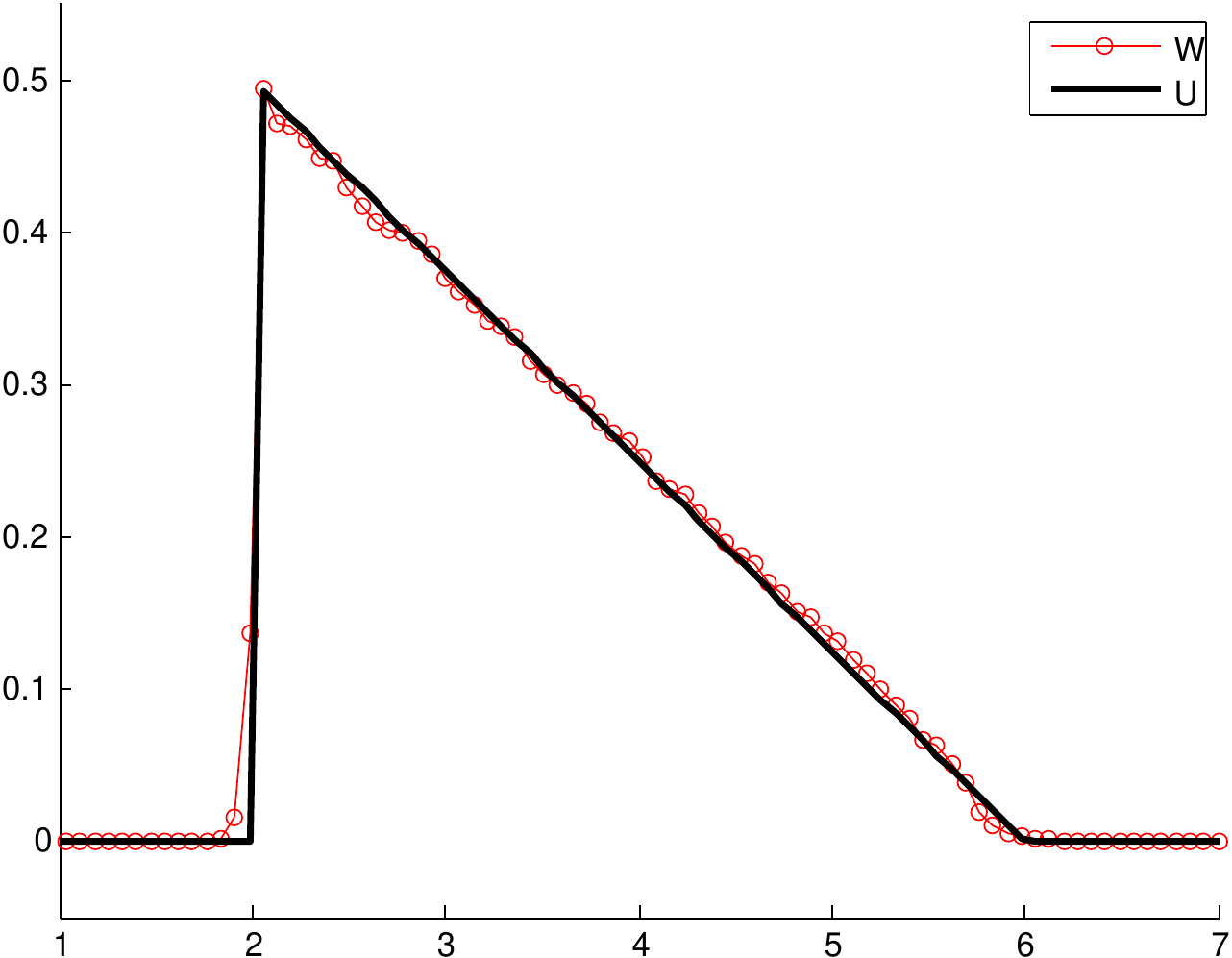} \\
(a) & (b)
\end{tabular}
\end{center}
\caption{Test \mMlwr: comparison with the exact solution. 
(a) $W^{N_T}(1,1)$ and $V^{N_T}(1,1)$. 
(b) $W^{N_T}(\thetauno^R,1)$.}
\label{fig:mMlwr:pures_and_bestR}
\end{figure}

In this test we solve equation \eqref{CL} with $f(u)=u(1-u)$ (LWR model for traffic flow). The exact solution is $u(x,T)=-\frac18 x+\frac34$ in $[2,6]$ and $u(x,T)=0$ elsewhere. A rarefaction and a shock are present.
The reference error (GODUNOV) in this case is equal to 0.0839. 
Richardson extrapolation method (with $s=\frac12$) suggests $\thetauno^R=0.956$ which leads to the error 0.0317 (-62.22\%). Here we used $W$ to evaluate the particles' velocity, see \eqref{fnl_WoV}.
Fig.\ \ref{fig:mMlwr:pures_and_bestR} shows the exact solution $U^{N_T}$ together with (a) the two uncoupled solutions $W^{N_T}(1,1)$ and $V^{N_T}(1,1)$, and (b) the best Richardson coupling $W^{N_T}(\thetauno^R,1)$ at final time.

%%%%%%%%%%%%%%%%%%%%%%%%%%%%%%%%%%%%%%%
%%%%%%%%%%%%%%%%%%%%%%%%%%%%%%%%%%%%%%%
%%%%%%%%%%%%%%%%%%%%%%%%%%%%%%%%%%%%%%%
\section*{Conclusions and perspectives}
We proposed a new method for coupling two (or more) explicit schemes approximating advection equations and scalar conservation laws. We provided some convergence analysis and we showed that effective couplings can be obtained in practice using a standard Richardson extrapolation. 
Moreover, several numerical tests confirmed that, whenever the constitutive schemes are respectively an Eulerian (macroscopic) and a Lagrangian (particle based) scheme, the blended scheme gives rather accurate results despite its Eulerian part gives only a low order approximation. 
This is a remarkable advantage, since it is well known that implementing high-order solvers is much easier for ODEs than for PDEs. 

It is interesting to note that the blended scheme can be put in relation with very different schemes proposed in the literature, see Section \ref{sec:intro}.
\begin{itemize}
\item If the coupling parameters are allowed to depend on space and/or on the solution itself, we can blend Eulerian schemes of different orders to recover the constitutive ideas of the \textit{filtered schemes}. 

\item Coupling an Eulerian scheme with a pure Lagrangian scheme and localizing particles, we recover a scheme that resembles that of the \textit{particle level-set method}.

\item Annealing the Eulerian scheme but keeping alive the coupling with $(\thetauno,\thetadue)=(0,1)$, we recover a scheme that resembles that of the \textit{particle-in-cell method}.

\item The multiscale blended scheme (with $\thetadue=1$) seems to follow the opposite philosophy with respect to the \textit{smoothed-particle hydrodynamics method}. While the latter regularizes the Lagrangian solution by a suitable choice of smoothing kernels, the former ``Dirac-izes'' the Eulerian solution. 
Moreover, the computation of the optimal blending parameter resembles the computation of the optimal supports for the smoothing kernels.
\end{itemize}

Finally, let us mention that one of the most promising and yet unexplored feature of the blended scheme is the \textit{parallelization on distributed-memory architectures}. Particle evolution is indeed embarrassingly parallel (within a time step), since there is no need to locate the neighbouring particles of a given one (which can be a severe bottleneck if particles live in different memory units).

%%%%%%%%%%%%%%%%%%%%%%%%%%%%%%%%%%%%%%%
%%%%%%%%%%%%%%%%%%%%%%%%%%%%%%%%%%%%%%%
%%%%%%%%%%%%%%%%%%%%%%%%%%%%%%%%%%%%%%%
\section*{Acknowledgements}
The authors thank R. Natalini who, commenting the paper \cite{cristiani2011MMS}, one day asked: ``Does the microscopic coupling decrease the numerical diffusion?''. E. Cristiani also wants to thank S. Succi for the motivating discussions and A. Di Mascio for the useful suggestions.

%%-----------------------------
%%      bibliography
%%-----------------------------

\end{document}